\documentclass{article}
\usepackage{dutchcal}

\newcommand{\mr}[1]{\mathring{#1}}
\usepackage{graphics}
  \usepackage{amssymb}
    \usepackage{amsthm}
  \usepackage{amsmath}
  \usepackage{pgf}
\usepackage{tikz}
\usetikzlibrary{arrows,automata}
\usepackage{circuitikz}
\usepackage{epstopdf}
  
\def\p{\partial}

\newtheorem{theorem}{Theorem}[section]
\newtheorem{lemma}[theorem]{Lemma}

\newtheorem{proposition}[theorem]{Proposition}
\theoremstyle{definition}
\newtheorem{definition}[theorem]{Definition}
\newtheorem{example}[theorem]{Example}
\newtheorem{remark}[theorem]{Remark}

\newcommand {\A}{\mathbb{A}}
\newcommand {\bay}{\begin{array}}
\newcommand {\eay}{\end{array}}
\newcommand {\bdm}{\begin{displaymath}}
\newcommand {\edm}{\end{displaymath}}
\newcommand{\ups}{\upsilon}
\newcommand{\w}{\varpi}

\newcommand{\vhi}{\varphi}

  \makeatletter
  \let\c@equation\undefined
  \let\c@section\undefined
  \let\c@subsection\undefined
  \let\c@zad\undefined
  \makeatother
  \newcounter{section}

  \newcounter{equation}[section]
  
  \newcounter{subsection}[section]

\newcommand{\beq}{\begin{equation}}
\newcommand{\eeq}{\end{equation}}

\newcommand{\N}{\mathbb N}

\newcommand{\bd}{\begin{displaymath}}
\newcommand{\ed}{\end{displaymath}}
\newcommand{\be}{\begin{equation}}
\newcommand{\ee}{\end{equation}}
\newcommand{\bdow}{\begin{proof}}
\newcommand{\edow}{\end{proof}}

\newcommand{\papap}{\papa \end{proof}}

\newfont{\smoldita}{cmmib8}
\newfont{\boldita}{cmmib10}
\newfont{\bboldita}{cmmib10}

\newcommand{\nn}{\nonumber}
\newcommand{\e}{\epsilon}

\newcommand {\sem}[1]{\mbox{$({#1}(t))_{t \geq 0}$}}

\newcommand{\cl}[2]{\int\limits_{#1}^{#2}}

\newcommand{\ti}[1]{\tilde{#1}}

\newcommand{\la}{\lambda}

\newcommand{\mb}[1]{\boldsymbol{#1}}
\newcommand{\mbb}[1]{\mathbb{#1}}
\newcommand{\mc}[1]{\mathcal{#1}}

\newcommand{\msf}[1]{\mathsf{#1}}

\newcommand{\bv}{\boldsymbol{v}}

\begin{document}
\begin{center}{\Large Telegraph systems on networks and port-Hamiltonians. I. Boundary conditions and well-posedness.}\end{center}

\begin{center}{J. Banasiak\footnote{The research has been partially supported by the National
Science Centre of Poland Grant 2017/25/B/ST1/00051 and the National Research Foundation of South Africa Grant 82770} \\ \small{Department of Mathematics and Applied Mathematics, University of Pretoria}\\ \small{Institute of Mathematics,  \L\'{o}d\'{z} University of Technology}\\ \small{International Scientific Laboratory of
Applied Semigroup Research, South Ural
State University}\\ \small{e-mail: jacek.banasiak@up.ac.za}\\\& \\A. B\l och\footnote{The research was completed while the author was a Doctoral Candidate in the Interdisciplinary Doctoral School at the \L\'{o}d\'{z} University of Technology, Poland.}\\
 \small{Institute of Mathematics,  \L\'{o}d\'{z} University of Technology} \\\small{e-mail: adam.bloch@dokt.p.lodz.pl}}\end{center}
\begin{abstract}
The paper is concerned with a system of linear hyperbolic differential equations on a network coupled through general transmission conditions of Kirchhoff's type at the nodes. We discuss the reduction of such a problem to a  system of 1-dimensional hyperbolic problems for the associated Riemann invariants and provide a semigroup theoretic proof of its well-posedness. A number of examples showing the relation of our results with recent research is also provided. 
\end{abstract}
\section{Introduction}

 In this paper we are concerned with dynamics described by a system of $2\times 2$ linear hyperbolic differential equations, defined on a collection of disconnected intervals and coupled through transmission conditions at the intervals' endpoints. A useful interpretation of such a system is a dynamical system on a metric graph, where the edges are identified with the intervals and thus the network problem is converted into a system of equations on an interval. We note that the restriction to $2\times 2$ systems is purely for convenience; the theory can be extended to hyperbolic systems of any (finite) dimension at each edge.

 Dynamical systems on graphs, or networks, have been studied since the early 80s, beginning with the pioneering paper \cite{Lum1}. They received more prominence in the context of quantum graphs, see \cite{Kuch, DM} and extensive references therein, where the dynamics on the edges are given by the heat or Schr\"{o}dinger equations. More general study of the diffusion on graphs can be found in \cite{vB, AB12, BFN3}.  Wave equations, both linear and nonlinear, on networks were considered in \cite{AM1, AM2}, and a comprehensive study of first order transport equations was undertaken in a series of papers such as \cite{KS2005, DKNR, BN, BFN3}. On the other hand, systems of hyperbolic first order equations have not received much attention until recently. Here we note the paper \cite{Car} on linearized blood flow,  the papers on the momentum operator \cite{Exn} and recent works on more general hyperbolic problems on graphs such as \cite{Nic, KMN}; the latter contains a comprehensive bibliography of the subject.

 Let us consider the system
  \begin{equation}
 \p_t\mb p + \mc M\p_x\mb{p} + \mc N \mb{p} =0, \quad t>0, 0<x<1,
 \label{sys10}
 \end{equation}
 where $\mb p = (p_1,p_2)^T$, $\mc M = (M_{lk})_{1\leq l,k\leq 2}$ and $\mc N = (N_{lk})_{1\leq l,k\leq 2}$ are real, possibly $x$-dependent, matrices, with $\mc M$ being strictly hyperbolic, that is, $\mc M$ has two distinct, real and nonzero eigenvalues. This implies the diagonalizability of $\mc M$, that is, the existence of an invertible linear operator $\mc F$ such that $\mc F^{-1}\mc M\mc F$ is a diagonal matrix. The components of $\mc F\mb p$ are known as the Riemann invariants of the system.  In this way, \eqref{sys10} with $\mc N=0$  becomes a system of  decoupled equations for the Riemann invariants. If we consider a collection of hyperbolic systems \eqref{sys10} on a network, then it must be complemented by boundary conditions which, in the simplest, local, case, express the relations between the values of the solutions along the edges at the nodes connecting these edges. The diagonalization of the differential operators on each edge leads to a system of decoupled differential equations, with boundary conditions for the resulting Riemann invariants coupled at the vertices of incident edges. Since the Riemann invariants have well-defined directions at which they move through the domain, we can talk about their incoming and outgoing values at a given endpoint. Given that each edge was parametrized and identified with the interval $(0,1)$, the boundary conditions express a relation between the incoming and outgoing values of the Riemann invariants (and thus of the values of the solution)   at $0$ and $1$. In this way, the problem becomes a special case of general 1-d hyperbolic  systems, \cite{BaCor}, or first order equations on a one dimensional spatial domain, \cite{Zwart2010, JaMoZw}, that recently have gained the name of (linear) first order port-Hamiltonian systems
 \cite[Definition 7.1.2]{JaZwbook}. Since in the linear case the differential equations of the system do not create much problem, the main issue in the analysis is ensuring that the general system of linear equations coupling the boundary values at 0 and 1 determines a sufficient (and necessary) number of the outgoing values to make the whole problem well-posed. We note that in, e.g. \cite{Diag}, the boundary conditions are written in an explicit form, already solved for the outgoing values at the endpoints. The general form of local linear boundary conditions is considered in \cite{Zwart2010, JaMoZw}, where, in contrast to \cite{Diag}, the authors also provide a fairly comprehensive solvability theory of such systems  based, however, on some advanced results from control theory, \cite[Chapter 7]{Staff}.

 One of the main remaining problems is to relate the results obtained for general port-Hamiltonians with the structure of the original network problem, where it is of importance to assign correct boundary conditions at the nodes of the network so that the resulting flow has some specified properties. Such conditions have appeared in e.g. \cite[Section 3]{Kuch1} in the context of quantum graphs, where the dynamics on each edge is given by the heat or Schr\"{o}dinger operators and it is required to construct boundary conditions ensuring the self-adjointness of the problem. Recently, similar ideas have been used in \cite{Nic, KMN}, where the authors considered hyperbolic systems at the edges of a network and at each node they constructed general boundary conditions that ensure the dissipativity of the problem in an appropriately weighted $L_2$-space. The construction is based on the integration by parts along each edge and ensuring that the end-point terms of the integration by parts are non-negative. The authors achieve this by assuming that the vertex values of the solution belong to the so-called totally isotropic subspace associated with the quadratic form defined by a symmetrization of the matrix $\mc M$. Due to this approach, the theory of \cite{Nic, KMN} is strongly dependent on the Hilbert space setting.

 Though the choice of the underlying state space for a problem is somewhat arbitrary, it should nevertheless be related to the physics of the problem as much as it is mathematically feasible. In fact, the same model can be analysed in different state spaces. Thus, to choose a proper mathematical setting, we should revisit the origins of the analysed model.

 \begin{example}
  System \eqref{sys10} is often referred to as the telegraph (or telegrapher's) system due to its interpretation as the system governing the electric voltage $p_1 =V$ and the current $p_2=I$ in a transmission line. In the simplest case of the lossless line there is neither resistance, nor inductance, hence $\mc N=0$, \cite[Section 7.7]{Fitz}. In this case, $M_{11} = M_{22} = 0$, $\frac{1}{M_{12}} = C$ and $\frac{1}{M_{21}}=L$ are, respectively,  the (constant) capacitance and the inductance per unit length of the line.  If we want to control the energy of the system, given by,
$$
E(t) = \frac{1}{2}\cl{0}{1} (CV^2(x,t) + LI^2(x,t))dx,
$$
\cite[p. 94]{JaZwbook}, then clearly the $L_2$ setting is physically meaningful. On the other hand, if we want to control the total charge $Q = CV$ and magnetic flux $\phi = LI,$ then the natural state space for the problem is $L_1$.
\end{example}

\begin{example}
\textbf{Shallow water waves -- Saint-Venant equations.} This is a system of nonlinear equations expressing the conservation of mass and momentum of a two-dimensional water flow in a horizontal channel of constant width. The system, derived in full generality in e.g. \cite[Appendix A]{LitFro},    can be written as
\begin{equation}
\begin{split}
\p_t \mc h   + \p_x (\mc u\mc h ) &=0,\\
\p_t(\mc u\mc h)+ \p_x\left(\mc u^2 \mc h + \frac{g}{2} \mc h^2\right) &=0,
\end{split}\label{SV1}
\end{equation}
where $\mc h$ is the depth of water and $\mc u$ its velocity at position $x$ and time $t$. The total mass $M$ and the total momentum $P$ of the water in a channel of length $L$ at a time $t$ are given by
$M(t) = \rho l\int_{0}^{L}\mc h(s,t)ds$ and $P(t)  = \rho l\cl{0}{L}\mc h(s,t)\mc u(s,t)ds$,
and thus \eqref{SV1} expresses the conservation of mass and momentum. It is more common in the literature to use $\mc h$ and $\mc u$ as the primary variables. Then, carrying out differentiations in brackets in \eqref{SV1} and simplifying, we get
\begin{equation}
\begin{split}
\p_t \mc h  + \mc h\p_x \mc u + \mc u\p_x\mc h &=0,\\
\p_t\mc u+  \mc u\p_x \mc u + g\p_x \mc h &=0.
\end{split}
\label{SV3}
\end{equation}
In many applications the Saint-Venant equations are linearized, typically at the  stationary state given, in this case, by a constant depth $H$ and velocity $V$ (and thus constant momentum $P = HV$). Considering small deviations from the equilibrium, $
u = \mc u - V, h = \mc h -H$ and $p =  \mc{uh}-P$, and ignoring higher order terms, we get
\begin{equation}
\begin{split}
\p_t  h  + H\p_x  u + V\p_x h &=0,\\
\p_t u+  V\p_x  u + g\p_x  h &=0.
\end{split}
\label{SV3l}
\end{equation}
The eigenvalues determining the directions of the flows are given by
\begin{equation}
\la_\pm = V\pm \sqrt{gH}.
\label{eigSV}
\end{equation}
The system is hyperbolic provided $gH\neq V^2$;  both eigenvalues are positive if the Froude number  $
Fr := {V}/{\sqrt{gH}} > 1$ and are of opposite sign if $Fr < 1$. In the latter case the flow is called subcritical, see \cite[Section 1.4]{BaCor}.

We observe that  the linearized momentum $p = \mc{uh} - P = (H+h)(V+u) - HV \approx Hu+Vh,$ up to higher order terms in $h$.
Thus the finite  $L_1(0,L)\times L_1(0,L)$ norm  of the solution describes the flow with finite mass and momentum  and thus $L_1(0,L)\times L_1(0,L)$ is a natural state space for \eqref{SV3l}.
\end{example}

\begin{example} \label{ex12}The system of equations describing a one-dimensional correlated random walk, \cite[Section 1.3.4]{BLbook} or \cite[Section 1.2]{Zaud}, takes the form
\begin{equation}
\begin{split}\p_t u&=-\gamma\p_x u-\lambda u+\lambda w,\\
\p_tw&=\gamma\p_x  w +\lambda u-\lambda w,
\end{split}
\label{w14}
\end{equation}
where $u$ and $w$ are the (probability) densities of particles moving, respectively, to the right and left, and $\la$ is the rate of the direction reversal. This system is often written in terms of  functions
\begin{equation} p(t,x)=u(t,x)+w(t,x),\qquad q(t,x)=u(t,x)-w(t,x),\label{w15}\end{equation}
which are the density and the net  current of the
particles. Adding and subtracting the equations in (\ref{w14}) we obtain
the telegrapher's system: \index{telegrapher's system}
\begin{equation}
\begin{split}
\p_tp+\gamma\p_x q&=0,\\
\p_tq+\gamma\p_x p +2\lambda q&=0.\end{split}\label{w16}\end{equation}\label{exrw}
Due to the probabilistic origins of the model, it is clear that the natural state space is the space of densities, a subset of $L_1(I)\times L_1(I)$, where $I\subset \mbb R$.
\end{example}
The presented examples show that it is important to define boundary conditions for hyperbolic network problems without relying on a particular state space and relate them with general boundary conditions of port-Hamiltonians and this is the aim of this paper and the companion paper, \cite{JBAB2}. Here, we define general Kirchhoff's type conditions at a vertex and provide conditions under which they determine all outgoing data from this vertex through the incoming ones, as required in \cite{BaCor}, and also that they  satisfy the solvability assumptions of \cite{JaMoZw}. We note that this contrasts with the approach of \cite{KMN}, where there is no a priori separation into the incoming and outgoing data in the boundary conditions at the cost of being confined to dissipative cases in the Hilbert space setting, not covering some standard cases and having not fully explicit representation, see Example \ref{KMNex}.     We show that our definition covers boundary conditions discussed in \cite{Nic, KMN} and can describe more general situations.  Next we provide an alternative, purely semigroup-theoretic, proof of the well-posedness of the problem. We start with the generation theorem in the $L_1$ setting, which is a straightforward adaptation of the approach in \cite{BN, BFN3}.  To extend the result to other $L_p$ spaces, $1<p<\infty,$ we use an  explicit representation of the constructed $L_1$ semigroup, valid for small times, and $L_p\subset L_1$, to show that the semigroup leaves $L_p$ invariant. In the paper \cite{JBAB2} we address the question when  port-Hamiltonian systems originate from hyperbolic network problems.
 \section{Notation and definitions}

 In this paper we shall be mainly concerned with real spaces. Let us take arbitrary $N\in \mbb N$. We consider $\mbb R^N$ with standard coordinate-wise partial order, $x = (x_1,\ldots,x_N)\geq y = (y_1,\ldots,y_N)$ if $x_i\geq y_i, i=1,\ldots,N.$ We write $x>y$ if $x\geq y$ and $x_i>y_i$ for at least one $i$ and $x\gg y$ if $x_i>y_i$ for any $i=1,\ldots,N$. For any $x\in \mbb R^N$ we denote $|x| = (|x_1|,\ldots,|x_N|)$ and, similarly, if $\mc B = (b_{ij})_{1\leq i\leq N, 1\leq j\leq M}$ is an $N\times M$ matrix, we write $|\mc B| = (|b_{ij}|)_{1\leq i\leq N, 1\leq j\leq M}.$

 We consider a network represented by a finite, connected and simple metric graph $\Gamma$ with $n$ vertices $\{\bv_j\}_{1\leq j\leq n}=:\Upsilon$ and $m$ edges $\{\mb e_j\}_{1\leq j\leq m}$. We denote by $E_{\bv}$ the set of edges incident to $\bv,$ let $J_{\bv} :=\{j;\; \mb e_j \in E_{\bv}\} $  and  $|E_{\bv}| =|J_{\bv}|$ be the valency of $\bv$. We identify the edges with unit intervals through sufficiently smooth invertible functions $l_j: \mb e_j \mapsto [0,1]$. In particular, we call $\bv$ with $l_j(\bv) =0$ the tail of $\mb e_j$ and the head if $l_j(\bv)=1$.  On each edge $\mb e_j$ we define $\mb p^j = (p^j_1, p^j_2)^T$ and consider a hyperbolic system
 \begin{equation}
 \p_t \mb p^j+ \mc M^j\p_x\mb p^j + \mc N^j \mb p^j =0, \quad t>0, 0<x<1, 1\leq j\leq m,
 \label{sys1}
 \end{equation}
 where $\mc M^j = (M^j_{lk})_{1\leq k,l\leq 2}$ and $\mc N^j = (N^j_{lk})_{1\leq k,l\leq 2}$ are real matrix functions defined on $[0,1]$, with $\mc M^j$ being strictly hyperbolic. More precisely, we assume that for each $x \in [0,1]$,
 $$
 \Delta^j = \left(\mathrm{tr \mc M^j}\right)^2 - 4\mathrm{det \mc M^j}>0
 $$
 so that the eigenvalues, given by
 $$
 \la^j_{\pm} = \frac{\mathrm{tr \mc M^j}\pm \sqrt{\left(\mathrm{tr \mc M^j}\right)^2 - 4\mathrm{det \mc M^j}}}{2},
 $$
 are real, different and nonzero, with $\la^j_-<\la^j_+$. We denote by $ f^j_\pm = (f^j_{\pm,1}, f^j_{\pm,2})^T$ the eigenvectors belonging to, respectively, $\la^j_\pm$.
 We observe that the eigenvalues can be of the same or of different signs. In the latter case, we have $\la^j_-<0<\la^j_+$.

 We introduce \begin{align*}
 \mc F^j &= \left(\begin{array}{cc} f^j_{+,1}&f^j_{-,1}\\
   f^j_{+,2}&f^j_{-,2} \end{array}\right),\qquad
 (\mc F^{j})^{-1} = \frac{1}{det \mc F^j}\left(\begin{array}{cc} f^j_{-,2}&-f^j_{-,1}\\
   -f^j_{+,2}&f^j_{+,1} \end{array}\right)
 \end{align*}
 and the Riemann invariants $\mb u^j = (u^j_1,u^j_2)^T, 1\leq j\leq m,$ by
 \begin{equation}
 \mb u^j = (\mc F^{j})^{-1} \mb p^j\quad\text{and}\quad
 \mb p^j = \binom{f^j_{+,1} u^j_1 + f^j_{-,1}u^j_2}{f^j_{+,2}u^j_1 + f^j_{-,2}u^j_2}.\label{pguw}
 \end{equation}
 Then we diagonalize \eqref{sys1} as
  \begin{align*}
 \p_t\mb u^j &= \p_t\left((\mc F^{j})^{-1} \mb{p}^j\right) = -(\mc F^{j})^{-1}\left(\mc M^j\p_x\mb{p}^j + \mc N^j \mb{p}^j\right)\nn\\
 &= \left(\begin{array}{cc}-\la^j_+&0\\0&-\la^j_-\end{array}\right)\p_x\mb{u}^j - (\mc F^{j})^{-1} \mc M^j (\p_x\mc F^{j}) \mb u^j - (\mc F^{j})^{-1}\mc N^j \mc F^{j}\mb u^j,
 \end{align*}
 for each $1\leq j\leq m$, or, in a compact form,
 \begin{equation}
 \p_t\mb u^j + \mc L^j \p_x\mb u^j + \overline{\mc N^j}\mb u^j= 0, \quad t>0, 0<x<1.
 \label{sysdiag}
 \end{equation}
Our assumptions ensure that $\mb u^j\mapsto \overline{\mc N^j}\mb u^j$ induces a bounded perturbation, hence all further considerations will be carried out for the case $\overline{\mc N^j}=0$.
\section{The boundary conditions}

Following the paradigm introduced in \cite[Section 1.1.5.1]{BaCor}, we require that at each point of the boundary all outgoing data must be determined by the incoming data. Since in a graph the boundary is represented by its vertices, the boundary conditions must represent some balance between the incoming and outgoing data at each vertex. However, for a general system  \eqref{sys1}, it is not always obvious which data are outgoing and which are incoming. For instance,  for the random walk model introduced in Example \ref{exrw}, it is easy to determine the outgoing and incoming data if the system is written in the form \eqref{w14} but it is not so clear for \eqref{w16} as $p$ is just the total density of particles and do not have any direction of the flow.  This can be resolved by using the Riemann invariants
$$
\mb u = \mc F^{-1} \mb p,
$$
where $\mb p = ((\mb p^j)_{1\leq j\leq m})^T,$  $\mb u = ((\mb u^j)_{1\leq j\leq m})^T$ and $
\mc F  = diag\{\mc F^j\}_{1\leq j \leq m}. $ If we disregard the lower order terms in \eqref{sysdiag}, equations for $\mb u_1 = (u^1_1,\ldots, u^m_1)^T$ and $\mb u_2 = (u^1_2,\ldots, u^m_2)^T$ are decoupled and, on $[0,1],$ the flow described by $u^j_1$ (respectively $u^j_2$)  occurs in the direction determined by the sign of $\la^j_+$ (respectively $\la^j_-$); that is, from 0 to 1 if the corresponding eigenvalue is positive, and from 1 to 0 otherwise. Let $\mb p \in (W^1_1(0,1))^{2m}$. Then, under the adopted assumptions on $\mc M^j, 1\leq j\leq m$,  also $\mb u \in (W^1_1(0,1))^{2m}$ and we can define the traces $\gamma_{\bv} \mb p^j$ and $\gamma_{\bv} \mb u^j$ whenever $\mb e_j \in E_{\bv}$. With some abuse of notation, we write $(\gamma_{\bv} \mb p^j)_{j\in J_{\bv}} = \mb p(\bv)$ for each $\bv\in \Upsilon$ and we use the same convention for $\mb u(\bv)$. We further allow $\bv$ to be replaced by $l_j(\bv)$  appropriate for  $\mb e_j$, $j \in J_{\bv}$. Since
$$
\mb p(\bv) = \mc F(\bv) \mb u(\bv),
$$
and  $\mc F(\bv)$ is invertible for any $\bv$, determining the boundary values  $\mb u(\bv)$ is equivalent to that for $\mb p(\bv)$  and thus we will mostly work with the former.
\begin{definition} Let $\mb v\in \Upsilon$. The outgoing values $u^j_k(\mb v), j\in J_{\bv}, k=1,2,$ are
 \begin{center}
 \begin{tabular} {|c|c|c|c|}
 \hline
 If &$\la^j_+>\la^j_->0$& $\la^j_+>0>\la^j_-$&$0>\la^j_+>\la^j_-$\\
 \hline
 $l_j(\bv) =0$& $u^j_1(\bv), u^j_2(\bv)$ & $u^j_1(\bv)$& none\\
 \hline
 $l_j(\bv) =1$&none &$u^j_2(\bv)$&$u^j_1(\bv), u^j_2(\bv)$ \\
 \hline
 \end{tabular}.
 \end{center}
  \label{tab1}
  \end{definition}
Denote by $\alpha_j$ the number of positive eigenvalues on $\mb e_j$. Then we see that for a given vertex $\bv$ with valence $|J_{\bv}|$ the number of outgoing values is given by
 \begin{equation}
 k_{\bv} := \sum\limits_{j\in J_{\bv}} (2(1-\alpha_j)l_j(\bv) +\alpha_j ).
 \label{kv1}
 \end{equation}

  \begin{definition} We say that $\bv$ is a \textit{sink}, and write $\bv\in \Upsilon_z$, if either $\alpha_j =2$ and $l_j(\bv) =1$ or $\alpha_j =0$ and $l_j(\bv) =0$ for all $j \in J_{\bv}$.  We say that $\bv$ is a \textit{source}, and write $\bv\in \Upsilon_s$, if either $\alpha_j =0$ and $l_j(\bv) =1$ or $\alpha_j =2$ and $l_j(\bv) =0$ for all $j \in J_{\bv}$. If $\bv$ is neither a source nor a sink, then we say that $\bv$ is a \textit{transient} (or \textit{internal}) vertex and write $\bv \in \Upsilon_t$.
\end{definition}
We observe that if $\bv\in \Upsilon_z$, then  $k_{\bv}=0$, while  if $\bv\in \Upsilon_s$, then  $k_{\bv}=2|J_{\bv}|$.
  Let us  introduce the partition  \begin{equation} \{1,\ldots,m\} = : J_1\cup J_2\cup J_0,\label{Jpart}\end{equation}
 where $j \in J_1$ if $\alpha_j =1, $  $j \in J_2$ if $\alpha_j =2$  and  $j\in J_0$ if $\alpha_j=0$. This partition induces the corresponding partition of each $J_{\bv}$ as
  $$J_{\bv} := J_{\bv,1}\cup J_{\bv,2}\cup J_{\bv,0}.
  $$
 We also consider another partition $J_{\bv} = J_{\bv}^0\cup J_{\bv}^1,$ where $j\in J_{\bv}^0$ if $l_j(\bv) =0$ and $j\in J_{\bv}^1$ if $l_j(\bv) =1$. Then we can give an alternative expression for $k_{\bv}$ as
 \begin{equation}
 k_{\bv} = \sum\limits_{j\in J^0_{\bv}} \alpha_j + \sum\limits_{j\in J^1_{\bv}} (2-\alpha_j) = |J_{\bv,1}|+ 2(|J^0_{\bv}\cap J_{\bv,2}| + |J^1_{\bv}\cap J_{\bv,0}|).
 \label{kval}
 \end{equation}
 The conditions of Definition \ref{tab1} can be alternatively expressed as follows.
   \begin{lemma} Let $\bv\in \Upsilon\setminus \Upsilon_z$. The outgoing boundary values of Definition \ref{tab1} are determined as follows
  \begin{description}
             \item (i) $u^j_1(0)$ is outgoing if and only if $j \in (J_{\bv,1}\cup J_{\bv,2})\cap J^0_{\bv},$
               \item (ii) $u^j_2(0)$ is outgoing if and only if $j \in J_{\bv,2}\cap J^0_{\bv},$
               \item (iii) $u^j_1(1)$ is outgoing if and only if $j \in J_{\bv,0}\cap J^1_{\bv},$
               \item (iv) $u^j_2(1)$ is outgoing if and only if $j \in (J_{\bv,1}\cup J_{\bv,0})\cap J^1_{\bv}.$
  \end{description}
  Hence, in particular,       \begin{equation}
    |(J_{\bv,1}\cup J_{\bv,2})\cap J^0_{\bv}| +|J_{\bv,2}\cap J^0_{\bv} | +
    |(J_{\bv,1}\cup J_{\bv,0})\cap J^1_{\bv}| +|J_{\bv,0}\cap J^1_{\bv}| = k_{\bv}.
    \label{kv3}
    \end{equation}\label{lemkv}
  \end{lemma}
  \begin{proof}
   By Definition \ref{tab1}, $u^j_1(0)$ is outgoing at $\bv$ if and only if $l_j(\bv) =0$ and either $\la^j_+>\la^j_->0$ or $\la^j_+>0>\la^j_-$ which occurs if and only if $j \in J_{\bv}^0$ and either $j\in J_{\bv,2}$ or $j \in J_{\bv,1}$ which is exactly (i). The  statements (ii)--(iv) follow in the same way. Eq. \eqref{kv3} follows from the one-to-one correspondence proved above.
  \end{proof}

\subsection{General Kirchhoff's conditions}
A typical example of the balance of incoming and outgoing data is Kirchhoff's law.  It has different interpretations, depending on the context. For an electrical circuit, it requires that the algebraic sum of currents at any node (vertex) must be zero, and thus expresses the charge conservation. On the other hand, for a flow in a channel network  it states that the rate of the fluid's inflow into any node must equal the rate of its outflow and thus  it is the mass conservation law.

Since Kirchhoff's law provides only one relation for the functions defined on the edges incident to $\mb v$,  typically it is not sufficient to specify all $k_{\mb v}$ outgoing values at $\mb v$,  which intuitively would lead to a well-posed initial boundary value problem for \eqref{sys1}.  Accordingly, for each $\bv \in \Upsilon\setminus\Upsilon_z$  we consider $k_{\bv}$ vectors of dimension $2|J_{\bv}|$ that, for a notational convenience, we write as
 $$
 \Phi_{\bv,r} = (\Phi_{\bv,r}^j)_{j \in J_{\bv}} =  ((\phi_{\bv,r}^j,\vhi_{\bv,r}^j))_{j \in J_{\bv}}, \quad r=1,\ldots,k_{\bv}.
 $$
We also introduce the matrix
\begin{equation}
\mb \Phi_{\bv}: = \left(\begin{array}{ccc}-&\Phi_{\bv,1} &-\\\vdots&\vdots&\vdots\\-&\Phi_{\bv,k_{\bv}}&-\end{array}\right) = \left(\begin{array}{ccccc}\phi^{j_1}_{\bv,1}&\vhi^{j_1}_{\bv,1}&\ldots&\phi^{j_{|J_{\bv}|}}_{\bv,1}&\vhi^{j_{|J_{\bv}|}}_{\bv,1}\\
\vdots&\vdots&\vdots&\vdots&\vdots\\
\phi^{j_1}_{\bv,k_{\bv}}&\vhi^{j_1}_{\bv,k_{\bv}}&\ldots&\phi^{j_{|J_{\bv}|}}_{\bv,k_{\bv}}&\vhi^{j_{|J_{\bv}|}}_{\bv,k_{\bv}}\end{array}\right),\
\label{Phiv}
\end{equation}
where $J_{\bv} = \{j_1, \ldots, j_{|J_{\bv}|}\}$.

\begin{definition}
We say that $\mb p$ satisfies a generalized Kirchhoff conditions at $\bv \in \Upsilon\setminus\Upsilon_z$ if
  \begin{equation}
  \sum\limits_{j \in J_{\bv}} (\phi^j_{\bv,r} p^j_1(\bv) + \vhi^j_{\bv,r} p^j_2(\bv))=0, \quad r=1,\ldots, k_{\bv},
  \label{bc2}
  \end{equation}
  or, in short,
  \begin{equation}
  \mb \Phi_{\bv} \mb p(\bv) = 0.
  \label{bc2s}
  \end{equation}
Equivalently,
\begin{equation}
\mb \Psi_{\bv} \mb u(\bv) : = \mb \Phi_{\bv} \mc F(\bv)\mb u(\bv)  = 0, \quad \bv \in \Upsilon\setminus\Upsilon_z.
  \label{bc2uw}
  \end{equation}
\end{definition}
  \subsubsection{Resolution at a vertex}
  To ensure that \eqref{bc2uw} satisfies the paradigm of \cite[Section 1.1.5.1]{BaCor} at each vertex, it must uniquely determine the outgoing values of $\mb u(\mb v)$ through the incoming values of $\mb u$ at $\mb v$.
Lemma \ref{lemkv} allows for an explicit formulation of the required assumptions. We introduce the block diagonal matrix
        \begin{equation}
    \mc {\ti {F}}_{out}(\bv) = diag\{\mc {\ti {F}}_{out}^j(\bv)\}_{j \in J_{\bv}},
    \label{frakVj}
    \end{equation}
    where
    $$
    \mc {\ti F}_{out}^j(\bv) = \left\{\begin{array} {ccc} \left(\begin{array}{cc}0&0\\0&0\end{array}\right)&\text{if} &j\in (J_{\bv,0} \cap J_{\bv}^0)\cup (J_{\bv,2}\cap J_{\bv}^1),\\
    \left(\begin{array}{cc}f^j_{+,1}(l_j(\bv))&f^j_{-,1}(l_j(\bv))\\f^j_{+,2}(l_j(\bv))&f^j_{-,2}(l_j(\bv))\end{array}\right)&\text{if}& j\in (J_{\bv,0} \cap J_{\bv}^1)\cup (J_{\bv,2}\cap J_{\bv}^0),\\
    \left(\begin{array}{cc}f^j_{+,1}(0)&0\\f^j_{+,2}(0)&0\end{array}\right)&\text{if} &j\in J_{\bv,1} \cap J_{\bv}^0,\\
    \left(\begin{array}{cc}0&f^j_{-,1}(1)\\0&f^j_{-,2}(1)\end{array}\right)&\text{if}& j\in J_{\bv,1} \cap J_{\bv}^1.
    \end{array}
    \right.
        $$
    Further, by $\mc F_{out}(\bv)$ we denote the contraction of $\mc {\ti F}_{out}(\bv)$; that is, the $2|J_{\bv}| \times k_{\bv}$ matrix obtained from $\mc {\ti F}_{out}(\bv)$ by deleting $2|J_{\bv}| - k_{\bv}$ zero columns, and then define $\mc F_{in}(\bv)$ as the analogous contraction of $\mc F(\bv)-\mc{\ti F}_{out}(\bv)$.

        In a similar way, we extract from $ \mb u(\bv)$ the outgoing boundary values $\widetilde{\mb u}_{out}(\bv) =(\widetilde{\mb u}^j_{out}(\bv))_{j\in J_{\mb v}},$ where
    \begin{align*}
        \widetilde{\mb u}^j_{out}(\bv) &= \left\{\begin{array} {ccc} (0,0)^T&\text{if} &j\in (J_{\bv,0} \cap J_{\bv}^0)\cup (J_{\bv,2}\cap J_{\bv}^1),\\
   (u^j_1(l_j(\bv)),u^j_2(l_j(\bv)))^T&\text{if}& j\in (J_{\bv,0} \cap J_{\bv}^1)\cup (J_{\bv,2}\cap J_{\bv}^0),\\
    (u^j_1(0),0)^T&\text{if} &j\in J_{\bv,1} \cap J_{\bv}^0,\\
    (0,u^j_2(1))^T&\text{if}& j\in J_{\bv,1} \cap J_{\bv}^1,
    \end{array}
    \right.
        \end{align*}
        and $        \widetilde{\mb u}_{in}(\bv) = {\mb u(\bv)}-\widetilde{\mb u}_{out}(\bv).
        $
As above, we define $\mb u_{out}(\bv)$ to be the vector in $\mbb R^{k_{\bv}}$ obtained by discarding the zero entries in $\widetilde{\mb u}_{out}(\bv)$ and, similarly, $\mb u_{in}(\bv)$ is the vector in $\mbb R^{2|J_{\bv}|-k_{\bv}}$ obtained from $\widetilde{\mb u}_{in}(\bv)$.
    \begin{proposition}
     Boundary system \eqref{bc2uw} at $\bv\in \Upsilon\setminus \Upsilon_z$ is equivalent to
    \begin{equation}
 \mb \Phi_{\bv} \mc F_{out}(\bv) \mb u_{out}(\bv)
    + \mb \Phi_{\bv} \mc F_{in}(\bv)\mb u_{in}(\bv) = 0
    \label{bcsplit1}
    \end{equation}
    and hence it  uniquely determines the outgoing values of $\mb u(\bv)$ at $\bv$ as defined by Definition \ref{tab1} if and only if
    \begin{equation}
    \mb\Phi_{\bv} \mc F_{out}(\bv)\quad \text{is\; nonsingular}.
    \label{compsolv}
    \end{equation}
     Then
    \begin{equation} \mb u_{out}(\bv) = - (\mb \Phi_{\bv} \mc F_{out}(\bv))^{-1} \mb \Phi_{\bv} \mc F_{in}(\bv)\mb u_{in}(\bv).
    \label{mcB}
    \end{equation}
     \label{prop1}
  \end{proposition}
  \begin{proof}  From  \eqref{frakVj} and the definition of $\widetilde{\mb u}_{out}(\bv)$ we see that
    \begin{equation}
    \mc{\ti F}_{out}(\bv)\widetilde{\mb u}_{in}(\bv) = 0, \qquad (\mc F(\bv)- \mc{\ti F}_{out}(\bv))\widetilde{\mb u}_{out}(\bv) = 0.
    \label{frakV}
    \end{equation}
          Hence
    $$
        \mc{\ti  F}_{out}(\bv){\mb u(\bv)} = \mc F_{out}(\bv)\mb u_{out}(\bv), \quad
     (\mc F(\bv)-\mc{\ti F}_{out}(\bv))\mb u(\bv) = \mc F_{in}(\bv)\mb u_{in}(\bv),
         $$
 and \eqref{bc2uw} can be written as
    \begin{align*}
    0&= \mb \Phi_{\bv} \mc F(\bv)\mb u(\bv)= \mb \Phi_{\bv} (\mc{\ti F}_{out}(\bv) + (\mc F(\bv) -\mc{\ti F}_{out}(\bv))(\widetilde{\mb u}_{out}(\bv)+\widetilde{\mb u}_{in}(\bv)) \\
    &= \mb \Phi_{\bv} \mc F_{out}(\bv) \mb u_{out}(\bv) + \mb \Phi_{\bv} \mc F_{in}(\bv)\mb u_{in}(\bv).
    \end{align*}
      \end{proof}

    \subsection{Graph independent boundary conditions}
    Assuming that vertices in $\Upsilon$ are ordered as $\{\bv_1, \ldots, \bv_n\}$, we can define $\mb \Psi' = diag \{\mb \Psi_{\bv} \}_{\bv\in \Upsilon\setminus\Upsilon_z},$  $\gamma \mb u  = ((\mb u(\bv))_{\bv \in \Upsilon\setminus\Upsilon_z})^T$ and write \eqref{bc2uw} in the global form
  \begin{equation}
  \mb \Psi' \gamma \mb u = 0.
  \label{bc2uwgl}
  \end{equation}
  By the hand shake lemma we have
  \begin{equation}
  2\sum\limits_{\bv\in \Upsilon} |J_{\bv}| = 4m
  \label{sumJv}
  \end{equation}
  and, by \eqref{kval} and $k_{\bv} =0$ for $\bv\in \Upsilon_z$,
  \begin{align*}
  \sum\limits_{\bv\in \Upsilon\setminus\Upsilon_z} k_{\bv} &= \sum\limits_{\bv\in \Upsilon} k_{\bv}=\sum_{\bv\in \Upsilon}|J_{\bv,1}|+ 2\left(\sum_{\bv\in \Upsilon}|J^0_{\bv}\cap J_{\bv,2}| + \sum_{\bv\in \Upsilon} |J^1_{\bv}\cap J_{\bv,0}|\right)\\& =\sum_{\bv\in \Upsilon}|J_{\bv,1}|+ \sum_{\bv\in \Upsilon}|J_{\bv,2}| + \sum_{\bv\in \Upsilon} |J_{\bv,0}|= 2m,
  \end{align*}
  where we used the fact that, when $\bv \in \Upsilon$ and $j \in J_{\bv}^0\cap J_{\bv,2}$ (so that $\bv$ is the tail of $\mb e_j$), there is exactly one $\bv'\in \Upsilon$ (the head of $\mb e_j$) such that $j \in J^1_{\bv'}\cap J_{\bv',2},$  so that
  $$
  2\sum_{\bv\in \Upsilon}|J^0_{\bv}\cap J_{\bv,2}| = \sum_{\bv\in \Upsilon}|J^0_{\bv}\cap J_{\bv,2}| + \sum_{\bv\in \Upsilon}|J^1_{\bv}\cap J_{\bv,2}| = \sum_{\bv\in \Upsilon}|J_{\bv,2}|;
  $$
   the same argument is valid for the second summand. Since \eqref{sumJv} contains sinks, it contains $2\sum\limits_{\bv\in \Upsilon_z} |J_{\bv}|$ irrelevant entries corresponding to the function values that are incoming at $\bv\in \Upsilon_z$ and do not influence any outgoing data. To keep, however, track of all vertex values, we augment  $\mb \Psi'$ to $\mb\Psi$ by adding zero columns corresponding to edges coming to sinks so $\mb\Psi$ is a $2m\times 4m$ matrix. In the same way, we can provide a global form of \eqref{bcsplit1}, splitting  \eqref{bc2uwgl} as
    \begin{equation}
 \mb \Psi^{out} \gamma \mb u_{out}
    + \mb \Psi^{in}\gamma\mb u_{in} = 0,
    \label{bcsplit2}
    \end{equation}
where $ \mb \Psi^{out}$ and  $\mb \Psi^{in}$ are $ diag \{\mb \Phi_{\bv} \mc F_{out}(\bv) \}_{\bv\in \Upsilon\setminus\Upsilon_z}$ and $diag \{\mb \Phi_{\bv} \mc F_{in}(\bv) \}_{\bv\in \Upsilon\setminus\Upsilon_z},$ augmented by zero columns corresponding to the incoming functions at the sinks, $\gamma \mb u_{out} = ((\mb u_{out}(\bv))_{\bv \in \Upsilon\setminus\Upsilon_z})^T$, and $\gamma \mb u_{in}$ is $((\mb u_{in}(\bv))_{\bv \in \Upsilon\setminus\Upsilon_z})^T$ augmented by incoming values at sinks.


  To give \eqref{bcsplit2} a vertex independent interpretation,  we focus on \eqref{sysdiag} with $\overline{\mc N^j}=0,$
  and discuss  general linear boundary conditions for $\mb u = ((u^j_1,u^j_2)_{1\leq j\leq m})^T$.
 Using the adopted parametrization, the flow described by $u^j_1$ occurs along $\mb e_j$ from the tail at $x=0$ to the head at $x=1$ if the corresponding eigenvalue, here $\la^j_+$,  is positive. Hence, using \eqref{Jpart}, if $j\in J_1$, then   we need to prescribe the value of $u^j_1$ at the tail, while $u^j_2$ there will be counted as incoming and the roles will be reversed at the head. Next, if $j\in J_2$, then both $u^j_1$ and $u^j_2$ flow from the tail to the head and both must be determined at the tail of the edge and provide the incoming information at its head; the picture will be reversed if $j \in J_0$.

This shows that we only need to distinguish functions describing the flow from $0$ to $1$ and from $1$ to $0$.
 Accordingly, we split $\mb u$ into parts corresponding to positive and negative eigenvalues:
 \begin{equation}
\begin{split}
\mb {\upsilon}  &:=\left((u^j_1)_{j\in J_1\cup J_2},(u^j_2)_{j\in J_2}\right) = (\ups_j)_{j\in J^+},\\
\mb {\w}& := \left((u^j_1)_{j\in J_0}, (u^j_2)_{j\in J_1\cup J_0}\right) = (\w_j)_{j\in J^-},
\end{split}
\label{renum}
\end{equation}
where, with some abuse of notation, $J^+ := \{J_1\cup J_2, J_2\}$ and $J^- := \{J_0, J_1\cup J_0\}$ are the sets of indices $j$ with, respectively, at least 1 positive eigenvalue, and at least 1 negative eigenvalue of $\mc M^j$. In  $J^+,$ respectively, $J^-$ the indices from $J_2$ (respectively $J_0$) appear twice so that we renumber them in some consistent way to avoid confusion. For instance, we can take  $J^+ =\{1,\ldots, m^u,m^u+1,\ldots, m^+\}$ and $J^- =\{m^++1,\ldots,m_u, m_u+1,\ldots, 2m\}$ and there are bijections between, respectively, $J_1\cup J_2$ and $\{1,\ldots, m^u\}$,
$J_2$ and $\{m^u+1,\ldots, m^+\}$, $J_0$ and $\{m^++1,\ldots,m_u\}$, and $J_1\cup J_0$ and $\{m_u+1,\ldots, 2m\}.$

We emphasize that different ways of indexing  would result in just re-labelling of the equations of \eqref{sysdiag} without changing its structure.

Note that in this way we converted the $2\times 2$ hyperbolic problem \eqref{sysdiag} on $\Gamma$ into a first order transport problem on a multi digraph $\mb \Gamma$ with the same vertices $\Upsilon$ and where each edge of $\Gamma$ was converted into two edges parametrized by $x\in [0,1],$ where $x=0$ and $x=1$ on both edges correspond  to the same vertices in $\Gamma$. Conversely, if we have a multigraph $\mb \Gamma,$ where all edges appear in pairs and  each two edges joining the same vertices are  parametrized concurrently, then we can collapse $\mb{\Gamma}$ to a graph $\Gamma$. Then \eqref{bcsplit2} can be written as
$$
\mb {\Xi} (\mb{\ups}(0), \mb{\ups}(1), \mb{\w}(0), \mb{\w}(1))^T = \mb{\Xi}_{out}(\mb{\ups}(0), \mb{\w}(1))^T + \mb{\Xi}_{in}(\mb{\ups}(1), \mb{\w}(0))^T= 0,
$$
where $\mb{ \Xi}$ is a $2m\times 4m$ matrix, and $\mb{\Xi}_{out}$ and $\mb{\Xi}_{in}$ are $2m\times 2m$ matrices, obtained as an appropriate permutation of columns of $\mb \Psi$ and $\mb{\Psi}^{out}$ and $\mb{\Psi}^{in}$, respectively.

We observe, however, that the above formulation  does not depend on the fact that $\mb{ \Xi}$ has a special form coming from local Kirchhoff's boundary conditions but, in principle, it can be an arbitrary $2m\times 4m$ matrix. Thus, we can consider
\begin{subequations}\label{ibvp0}
\begin{equation}
\p_t \mb{\ups}(x,t) = -c_j(x)\p_x \mb{\ups}(x,t),\quad t>0,0<x<1,\label{ibvpa0}
\end{equation}
\begin{equation}
\p_t \mb{\w}(x,t) = c_j(x)\p_x \mb{\w}(x,t), \quad t>0,0<x<1,\label{ibvpb0}
\end{equation}
\begin{equation}
\begin{split}
\mb{\ups}(x,0) &= \mr{\mb{\ups}}(x),\;0<x<1,\\
\mb{\w}(x,0) &= \mr{\mb{\w}}(x),\;0<x<1,  \label{ibvpc0}
\end{split}
\end{equation}
\begin{equation}
\mb {\Xi} (\mb{\ups}(0), \mb{\ups}(1), \mb{\w}(0), \mb{\w}(1))^T = 0, \quad t>0, \label{ibvpd0}
\end{equation}
\end{subequations}
where $\mb{\ups}(x,t) =(\ups_j(x,t))_{j\in J^+},$ $ \mb{\w}(x,t)= (\w_j(x,t))_{j\in J^-}$ and
  the positive functions $c^j, j\in J^+\cup J^-,$  equal the absolute values of the corresponding eigenvalues.  
\section{The generation theorem}
We observe that our network problem has become a first order problem and, in fact, a special case of the so-called port-Hamiltonian systems, see e.g. \cite{JaZwbook}. Though the definition in \cite{JaZwbook} concerns the Hilbert space case, it follows from \cite{JaMoZw} that the Hilbert space structure is needed to reduce a general port-Hamiltonian to the diagonal case \eqref{ibvpa0}, \eqref{ibvpb0}, while the well-posedness theory of the latter, developed in  \cite{Zwart2010, JaMoZw, EngKra}, applies in any $L_p$ space, $1\leq p <\infty$. This theory, however, is based on control theory results developed in \cite{Staff}. Here we present an alternative, semigroup theoretic proof, by placing \eqref{ibvp0}  in the framework of \cite{BFN3} that immediately leads to the well-posedness of \eqref{ibvp0} in the $L_1$ setting. Then the $L_p$ theory follows by direct estimates of the $L_1$ solutions with $L_p$ data. We observe that by further re-parametrizing of the arcs of $\mb \Gamma$ with $j\in J^-$ we could transform \eqref{ibvpa0}, \eqref{ibvpb0} to the case with positive transport speeds, making thus the calculations of \cite{BFN3} directly available. We decided, however, to leave the problem in the form with the separated directions of transport due to its natural connection with the original second order problem, and also since this is the form studied in \cite{BaCor, Zwart2010, JaMoZw}.

Hence, we consider an arbitrary $2m\times 4m$ real matrix $\mb \Xi$ in \eqref{ibvpd0} and split it into $2m\times 2m$ matrices as $
\mb \Xi = (\mb \Xi_{out}, \mb \Xi_{in})$, so that \eqref{ibvpd0} takes the form
\begin{equation}
\mb \Xi_{out}(\mb{\ups}(0,t),\mb{\w}(1,t)) =- \mb \Xi_{in}(\mb{\ups}(1,t),\mb{\w}(0,t)).
\label{bsys}
\end{equation}
By combining \cite[Theorem 1.5]{JaMoZw} with \cite[Theorem 3.3]{Zwart2010}, we see that proving the well-posedness of \eqref{ibvp0}, there is no loss of generality assuming that
\begin{equation}
\mb\Xi_{out}\quad \text{is\; invertible}
\label{Xiinv}
\end{equation}
as the existence of $\mb\Xi^{-1}_{out}$ is necessary for the semigroup generation.
\begin{remark} Clearly, if \eqref{compsolv} satisfied, then $\mb\Xi^{-1}_{out}$ exists. However, as noticed in \cite[Theorem 3.7]{KMN} (for dissipative boundary conditions in a Hilbert space setting), the solvability of $\mb\Psi^{out}_{\bv}:=\mb \Phi_{\bv} \mc F_{out}(\bv)$ at each $\mb v\in \Upsilon\setminus\Upsilon_z$ is not necessary as long as there is some global solvability and  \eqref{Xiinv} renders such a condition.
\end{remark}

Let us denote
\bd
\mc{C}(x):=diag\left\{\left(-c_j(x)\right)_{j\in J^+},\left(c_j(x)\right)_{j\in J^-}\right\} = diag \{-\mc C_+(x),\mc C_-(x)\}
\ed
and, for an arbitrary  $2m\times 4m$ real matrix $\mb\Xi$ satisfying \eqref{Xiinv}, let  $\mc B := \mb\Xi_{out}^{-1}\mb \Xi_{in}$. Let $\mathcal{B}=\left(B^{kl}\right)_{k,l=1,2}$, where $B^{kl}:=\left(b^{kl}_{ij}\right)_{i,j}$ with $i\in J^+$ if $k=1$ and $i\in J^-$ if $k=2$ and $j\in J^+$ if $l=1$ and $j\in J^-$ if $l=2.$ Now consider the spaces $\textbf{X}_p:=(L^p(0,1))^{2m}$ and $\mb Y_p = (W^p_1(0,1))^{2m}$ for $1\leq p<\infty,$ and define an operator $({A}_\mc B, D_p({A}_\mc B))$ in $\textbf{X}_p$ as $\msf A|_{D_p({A}_\mc B)}$, where
\begin{subequations}\label{opAB}
\begin{equation}
\mathsf{A}:=\mc{C}diag \{\underbrace{\partial_x,\ldots,\partial_x}_{2m\, \text{times}}\},
\label{opABa}
\end{equation}
\begin{equation}
D_p({A}_\mc B)\!:=\!\left\{\!\binom{\mb{\ups}}{\mb{\w}}\!\!\in \mb Y_p;\;\!\binom{\mb{\ups}(0)}{\mb{\w}(1)}\! =\! \mc B \binom{\mb{\ups}(1)}{\mb{\w}(0)}\right\}.
\label{opABb}
\end{equation}
\end{subequations}
Then we have the following theorem, whose proof closely follows that of \cite[Theorem 3.1]{BFN1} and thus will be only sketched.
\begin{theorem}
Let $\mc B$ be an arbitrary matrix. The operator $({A}_\mc B,D_1({A}_\mc B))$ generates a $C_0$-semigroup on $\mb X_1$.
\end{theorem}
\begin{proof}
Since $(C^{\infty}_0(0,1))^{2m}\subset D_1({A}_\mc B),$  $D_1({A}_\mc B)$ is dense in $\textbf{X}_1$.

Let $\mb f  \in\textbf{X}_1$. The resolvent equation for $A_\mc B$ takes the form
\begin{subequations}\label{reseq1}
\begin{equation}
\lambda \ups_j(x)+c_j(x)\partial_x\ups_j(x)=f_j(x), \quad j\in J^+\label{reseq1a}
\end{equation}
\begin{equation}
\lambda \w_j(x)-c_j(x)\partial_x\w_j(x)=f_j(x), \quad j\in J^-.\label{reseq1b}
\end{equation}
\end{subequations}
 Let us denote
 $$
 e_\la(a,b):=e^{-\lambda\int\limits_a^b\frac{1}{c_j(z)}\,dz},
 $$
 and $
\mc{E}_{\lambda}^{+}(a,b):=diag\Big\{e_\la(a,b)\Big\}_{j\in J^+}$, $\mc{E}_{\lambda}^{-}(a,b):=diag\Big\{e_\la(a,b)
\Big\}_{j\in J^-}.
$
For  $\mb \xi = (\xi_1,\ldots,\xi_{2m})$, we set $(\xi_j)_{j\in J^+} =\mb\xi^+$ and $(\xi_j)_{j\in J^-} =\mb\xi^-$. Then, integrating \eqref{reseq1a} from 0 to $x$ and \eqref{reseq1b} from $x$ to 1, the solution to \eqref{reseq1}  can be written as
\begin{equation}\label{solutiontotheresolventequation}
\begin{split}
\left(\begin{array}{cccc}
\mb\ups^+(x)\\
\mb\w^-(x)
\end{array}\right)&=diag \{\mc{E}_{\lambda}^{+}(0,x),  \mc{E}_{\lambda}^{-}(x,1)\} \left(\begin{array}{c}
\mb\ups^{0}\\
\mb\w^{0}
\end{array}\right)\\
&\phantom{xx}+\left(\begin{array}{c}
\int\limits_0^x\mc{E}_{\lambda}^{+}(s,x)\mc{C}_+^{-1}(s)\mb f^{+}(s)\,ds\\
\int\limits_x^1\mc{E}_{\lambda}^{-}(x,s)\mc{C}_-^{-1}(s)\mb f^{-}(s)\,ds\end{array}\right),
\end{split}
\end{equation}
for an arbitrary vector $(\mb\ups^0,\mb\w^0)^T.$
Using the boundary condition in $D_1({A}_\mc B)$ to determine $\mb\ups^{0},\mb\w^{0},$ we have
\begin{equation}
(\mc I-\mathcal{B}\mathcal{E}_{\lambda}(0,1))\left(\begin{array}{c}
\mb\ups^{0}\\
\mb\w^{0}
\end{array}\right)=\mathcal{B}\left(\begin{array}{c}
\int\limits_0^1\mc{E}_{\lambda}^{+}(s,1)\mc{C}_+^{-1}(s)\mb f^{+}(s)\,ds\\
\int\limits_0^1\mc{E}_{\lambda}^{-}(0,s)\mc{C}_-^{-1}(s)\mb f^{-}(s)\,ds
\end{array}\right),\label{ups}
\end{equation}
where $\mc I$ is the $2m\times 2m$ identity matrix and $\mathcal{E}_{\lambda}(0,1):=diag\left\{\mc{E}_{\lambda}^{+}(0,1),\mc{E}_{\lambda}^{-}(0,1)\right\}.$
If $M:=\min\limits_{j \in J^+\cup J^-} \left(\left\|{c_j^{-1}}\right\|_{L^1(0,1)}\right)
$, then $\left\|\mathcal{B}\mathcal{E}_{\lambda}(0,1)\right\|<1$ for any ${\lambda>\frac{\ln\left\|\mathcal{B}\right\|}{M}},$
hence
\be\label{neumannseriesrepresentation}
\left(\begin{array}{c}
\mb\ups^{0}\\
\mb\w^{0}
\end{array}\right)=\sum\limits_{n=0}^{\infty}\left(\mathcal{B}\mathcal{E}_{\lambda}(0,1)\right)^n\mathcal{B}\left(\begin{array}{c}
\int\limits_0^1\mc{E}_{\lambda}^{+}(s,1)\mc{C}_+^{-1}(s)\mb f^{+}(s)\,ds\\
\int\limits_0^1\mc{E}_{\lambda}^{-}(0,s)\mc{C}_-^{-1}(s)\mb f^{-}(s)\,ds
\end{array}\right),
\ee
so the resolvent  $R(\lambda,A_\mc B)$ is fully determined by \eqref{solutiontotheresolventequation} and \eqref{neumannseriesrepresentation}.

Let us define a modified operator ${A}_{\left|\mathcal{B}\right|}$ as in \eqref{opAB} but with $\mc B$ replaced by $|\mc B|$. Then, by \eqref{solutiontotheresolventequation} and \eqref{neumannseriesrepresentation},   $R(\lambda,{A}_{\left|\mathcal{B}\right|})$ is a positive operator such that
$$
\left|R(\lambda,{A}_\mc B)
\mb f\right|\leq R(\lambda,{A}_{\left|\mathcal{B}\right|})\left|\mb f\right|.
$$
Iterating this inequality and using the lattice property of the norm,  we have
\begin{equation}
 \left\|R(\lambda,{A}_\mc B)^n\mb f\right\|\leq \left\|R(\lambda,{A}_{\left|\mathcal{B}\right|})^n\right\|\left\|\mb f\right\|, \quad  n \geq 1. \label{resineq}
\end{equation}
Hence if ${A}_{\left|\mathcal{B}\right|}$ generates a semigroup, so does ${A}_\mc B$ and so, from now on, we shall consider operators ${A}_{\mathcal{B}}$ with $\mc B\geq 0$ so that $R(\la, A_{\mc B})$ is  positive and thus we can use $\mb f \geq 0$ for the norm estimates.

Let us denote
\be\label{columnsums}
b_j^{kl}:=\sum\limits_{i}b_{ij}^{kl}, \ \ k,l=1,2,
\ee
where $i$ runs through $J^+$ if $k=1$, and it runs through $J^-$ if $k=2$. That is,  $b_j^{kl}$ is the sum of elements in the $j$-th column of the matrix $B^{kl}$.
Adding separately the equations with indices in $J^+$ and $J^-$ in \eqref{ups}  and using \eqref{columnsums}, we have
\begin{equation}\label{sumofuzeroplus}
\begin{split}
&\sum\limits_{i\in J^+}\ups^0_i=\!\!\sum\limits_{j\in J^+}b_{j}^{11}e_\la(0,1)\ups^0_j+\sum\limits_{j\in J^+}b_{j}^{11}\int\limits_0^1\frac{e_\la(s,1)}{c_j(s)}f_j(s)\,ds\\
&\phantom{xxxxxxx}+\sum\limits_{j\in J^-}b_{j}^{12}e_\la(0,1)\w^0_j+\sum\limits_{j\in J^-}b_{j}^{12}\int\limits_0^1\frac{e_\la(0,s)}{c_j(s)}f_j(s)\,ds
\end{split}
\end{equation}
and
\begin{equation}\label{sumofuzeroplus1}
\begin{split}
&\sum\limits_{i\in J^-}\w^0_i=\!\!\sum\limits_{j\in J^+}b_{j}^{21}e_\la(0,1)\ups^0_j+\sum\limits_{j\in J^+}b_{j}^{21}\int\limits_0^1\frac{e_\la(s,1)}{c_j(s)}f_j(s)\,ds\\
&\phantom{xxxxxxx}+\sum\limits_{j\in J^-}b_{j}^{22}e_\la(0,1)\w^0_j+\sum\limits_{j\in J^-}b_{j}^{22}\int\limits_0^1\frac{e_\la(0,s)}{c_j(s)}f_j(s)\,ds.
\end{split}
\end{equation}
Now let us introduce a new, equivalent, norm in $\textbf{X}_1$  by
\bd
\left\|\left(\mb\ups,\mb\w\right)\right\|_c:=\sum\limits_{j\in J^+}\left\|{\ups_j}{c^{-1}_j}\right\|_{L^1(0,1)}+\sum\limits_{j\in J^-}\left\|{\w_j}{c^{-1}_j}\right\|_{L^1(0,1)}.
\ed
With this new norm, for $\mb f\geq 0$, we have by \eqref{solutiontotheresolventequation},
\begin{align*}
\left\|R(\lambda,{A}_{\mathcal{B}})\mb{f}\right\|_c&=\sum\limits_{j\in J^+}\int\limits_0^1\frac{e_\la(0,x)}{c_j(x)}\ups^0_j\,dx
+\sum\limits_{j\in J^+}\int\limits_0^1\int\limits_0^x\frac{e_\la(s,x)}{c_j(x)}\frac{f_j(s)}{c_j(s)}\,ds\,dx\\
&\phantom{x}+\sum\limits_{j\in J^-}\int\limits_0^1\frac{e_\la(x,1)}{c_j(x)}\w^0_j\,dx
+\sum\limits_{j\in J^-}\int\limits_0^1\int\limits_x^1\frac{e_\la(x,s)}{c_j(x)}\frac{f_j(s)}{c_j(s)}\,ds\,dx\\
&=
\frac{1}{\lambda}\sum\limits_{j\in J^+}\left(\ups_j^{0}-\ups_j^{0}e_\la(0,1)+\int\limits_0^1\frac{f_j(s)}{c_j(s)}\,ds
-\int\limits_0^1\frac{e_\la(s,1)}{c_j(s)}f_j(s)\,ds\right)\\
&+\frac{1}{\lambda}\sum\limits_{j\in J^-}\left(\w^0_j
-\w_j^0e_\la(0,1)+\int\limits_0^1\frac{f_j(s)}{c_j(s)}\,ds
-\int\limits_0^1\frac{e_\la(0,s)}{c_j(s)}f_j(s)\,ds\right),
\end{align*}
where we changed the order of integration and used
\bd
\frac{d}{dx}e^{-\lambda\int\limits_0^x\frac{1}{c_j(z)}\,dz}=-\frac{\lambda}{c_j(x)}e^{-\lambda\int\limits_0^x\frac{1}{c_j(z)}\,dz}
\ed
to evaluate the resulting inner integrals. Now we make use of \eqref{sumofuzeroplus} and \eqref{sumofuzeroplus1} to obtain, after some algebra,
\begin{align*}
&\left\|R(\lambda,{A}_{\mathcal{B}})\mb f\right\|_c=\frac{1}{\lambda}\sum\limits_{j\in J^+}\left(b_j^{11}+b_j^{21}-1\right)\left(\ups_j^{0}e_\la(0,1)+\int\limits_0^1\frac{e_\la(s,1)}{c_j(s)}f_j(s)\,ds\right)\\
&\phantom{xxxxxx}+\frac{1}{\lambda}\sum\limits_{j\in J^-}\left(b_j^{12}+b_j^{22}-1\right)\left(\w_j^{0}e_\la(0,1)+\int\limits_0^1\frac{e_\la(0,s)}{c_j(s)}f_j(s)\,ds\right)+\frac{1}{\lambda}\left\|
\mb f\right\|_c.
\end{align*}
Observe that only the expressions involving $b_j^{kl}, k,l=1,2$ can be negative. Thus we consider three cases.\\
\textit{Case 1.} $b_j^{11}+b_j^{21}\leq 1, b_j^{12}+b_j^{22}\leq 1$ for all $j\in J^+\cup J^-$.
Then we have
	\bd
	\left\|R(\lambda,{A}_{\mathcal{B}})\mb f \right\|_c\leq \frac{1}{\lambda}\left\|\mb f\right\|_c
	\ed
	and hence ${A}_{\mathcal{B}}$ generates a positive semigroup of contractions in $\left(\textbf{X}_1,\left\|\cdot\right\|_c\right)$ and, by the equivalence of norms, a positive bounded semigroup in $\textbf{X}_1$.\\
\textit{Case 2.} $b_j^{11}+b_j^{21} \geq 1, b_j^{12}+b_j^{22}\geq  1$ for all $j\in J^+\cup J^-$. Then
	\bd
	\left\|R(\lambda,{A}_{\mathcal{B}})\mb f\right\|_c\geq \frac{1}{\lambda}\left\|\mb
f\right\|_c
	\ed
	and, by \cite[Theorem 2.5]{Arendt},  ${A}_{\mathcal{B}}$ generates a~positive semigroup in $\left(\textbf{X}_1,\left\|\cdot\right\|_c\right),$ and hence in $\textbf{X}_1$.\\
	 \textit{Case 3.} There exist subsets $I'\subset J^+$ and $I''\subset J^-,$ such that at least one is nonempty, $I'\cup I'' \neq J^+\cup J^-$ and
	\begin{equation}\label{condi}
\begin{split}
	 b^{11}_j+b^{21}_j<1, \quad j\in I', &\qquad  b^{11}_j+b^{21}_j\geq 1,\quad j\in J^+\setminus I',\\
		b^{12}_j+b^{22}_j<1, \quad j\in I'', &\qquad  b^{12}_j+b^{22}_j\geq 1, \quad j\in J^-\setminus I''.
	\end{split}
\end{equation}
	We introduce a new matrix $\tilde{\mathcal{B}}$ by replacing the entries giving rise to $b^{11}_j+b^{21}_j<1, j\in I',$ and $b^{12}_j+b^{22}_j<1, j\in I'',$ by 1 and leaving the other entries unchanged.  Then $\tilde{\mathcal{B}}$ satisfies the assumption of Case 2 and $\mathcal{B}\leq\tilde{\mathcal{B}}$ as \eqref{condi} implies that each replaced entry of $\mc B$ was smaller than 1 (and positive). As before, we introduce the operator ${A}_{\tilde{\mathcal{B}}}$ on $D_1(A_{\mc{\ti B}})$.
	Repeating the argument leading to \eqref{resineq}, there exists $R(\lambda,{A}_{\tilde{\mathcal{B}}})$ for sufficiently large $\lambda$ and
	$	0\leq R(\lambda,{A}_{\mathcal{B}})\leq R(\lambda,{A}_{\tilde{\mathcal{B}}}),$  hence
			\bd
	\left\|R(\lambda,{A}_{\mathcal{B}})^n\right\|\leq \left\|R(\lambda,{A}_{\tilde{\mathcal{B}}})^n\right\|, \quad {n\in\N},
	\ed
so also ${A}_{\mathcal{B}}$ is the generator.
	
We have shown that for any nonnegative matrix $\mc B$, the operator ${A}_{\mathcal{B}}$ generates a~semigroup. Therefore, by \eqref{resineq},   this is  true for ${A}_\mc B$ with arbitrary matrix $\mc B$.
\end{proof}
Next we  show that we can use the above result to prove the generation in any $\mb X_p$ with $p\in [1,\infty)$.
Let \sem{G_1} be the semigroup on $\mb X_1$ generated by $A_{1,\mc B}$.
\begin{theorem} The restriction $\sem{G_p}=(G_1(t)|_{\mb X_p})_{t\geq 0}$ is a strongly continuous semigroup on $\mb X_p$ whose generator is the part of $A_{1,\mc B}$ in $\mb X_p$.
\end{theorem}
\begin{proof}
If $(\mr {\mb\ups},\mr{\mb\w})\in D(A_{1,\mc B})$, then  $(\mb\ups(x,t), \mb\w(x,t)) = [G_1(t)(\mr {\mb\ups},\mr{\mb\w})](x)$ is an absolutely continuous solution to \eqref{ibvp0}.
 Let us consider $j\in J^+$ and denote  $X_j(x,t) = L_{j}^{-1}(L_j(x) -t),$ where
 \begin{equation}
 L_j(x) = \cl{0}{x}\frac{ds}{c_j(s)}.
 \label{Lj}
 \end{equation}
 We have $L_j:[0,1]\mapsto [0, L_j(1)]$ (where $L_j(1) =: T_j$ is the  time needed to traverse the edge $j$ with the speed $c_j(x)$ from the tail at $x=0$ to the head at $x=1$); it is a strictly increasing function and hence its inverse $L_{j}^{-1}: [0, T_j]\mapsto [0,1]$ is well defined. Then, by the uniqueness of solutions to scalar first order partial differential equations,
 $$
 \ups_j(x,t) = \mr {\ups}_j(X_j(x,t)), \quad 0\leq L_j(x)-t\leq T_j.
 $$
 In particular, $\ups_j(1,t) = \mr {\ups}_j(X_j(1,t)) = \mr {\ups}_j(L_j^{-1}(T_j-t))$ is defined for $0\leq t\leq T_j.$

 Similarly, for $j\in J^-$
 $$
  \w_j(x,t) = \mr {\w}_j(X_j(x,-t)), \quad 0\leq L_j(x)+t\leq T_j
 $$
 and $\w_j(0,t) = \mr {\w}_j(X_j(0,-t)) = \mr {\w}_j(L^{-1}_j(t)) $ is defined for  $0\leq t\leq T_j.$

 Let us define $X_{k,j}(x,t) = L^{-1}_k(L_j(x)-t)$ and $Y_{k,j}(x,t) = L^{-1}_k(t-L_j(x))$, $T= \min_{j \in J^+\cup J^-}\{ T_j\}$. We fix $j\in J^+.$  Thus, by the uniqueness of solutions, for $t\in [0,T]$ and $(\mr{\mb\ups}, \mr {\mb\w}) \in D(A_{1,\mc B})$ we have
  \begin{equation}
 \ups_j(x,t) = \left\{\begin{array}{lcl} \mr {\ups}_j(X_j(x,t)), && x\in (L^{-1}_j(t),  1], \\
 \sum\limits_{k \in J^+}b_{jk}\mr{\ups}_k (X_{k,j}(x,t-T_k))&&\\
  + \sum\limits_{k\in J^-}b_{jk}\mr{\w}_k(Y_{k,j}(x,t)), && x\in [0, L^{-1}_j(t)),
   \end{array}
 \right.
 \label{upsdef}
 \end{equation}
 and, similarly for $j\in J^-$,
  \begin{equation}
 \w_j(x,t)\! =\! \left\{\!\!\begin{array}{lcl} \mr {\w}_j(X_j(x,-t)), &&\!\!  x\in [0, L_j^{-1}(T_j-t)), \\
 \sum\limits_{k \in J^+}b_{jk}\mr{\ups}_k (Y_{k,j}(x, T_k +T_j-t)))&&\\
  + \sum\limits_{k\in J^-}b_{jk}\mr{\w}_k(X_{k,j}(x, T_j-t)), && \!\!x\in (L_j^{-1}(T_j-t),1].
 \end{array}
 \right.
 \label{wdef}
 \end{equation}
 Let us denote by $(\mr{\mb\ups},\mr{\mb \w}) \mapsto \mc G(t)(\mr{\mb \ups},\mr{\mb \w}), 0\leq t\leq T,$ the operator defined by \eqref{upsdef} and \eqref{wdef}. Since the functions $L_j, j \in J^+\cup J^-,$ are diffeomorphisms, the compositions  with functions $(\mr\ups, \mr \w) \in \mb X_p$  are measurable.
  Then, for $(\mr\ups, \mr \w) \in \mb X_p, 0\leq t\leq T$ and $j \in J^+$, we have
    \begin{align}\label{upse1}
    &\cl{0}{1} |\ups_j(x,t)|^pdx \leq \cl{L_j^{-1}(t)}{1}|\mr {\ups}_j(X_j(x,t))|^pdx\nn\\& + D^{p/q}_b \left( \sum\limits_{k \in J^+}\cl{0}{L_j^{-1}(t)} |\mr{\ups}_k (X_{k,j}(x,t-T_k))|^pdx +\sum\limits_{k \in J^-}\cl{0}{L_j^{-1}(t)}|\mr{\w}_k(Y_{k,j}(x,t))|^pdx \right)\nn\\
    &= \cl{0}{L_j^{-1}(T_j-t)}|\mr {\ups}_j(z)|^p\frac{c_j(X_j^{-1}(z,t))}{c_j(z)}dz\\& + D^{p/q}_b \left( \sum\limits_{k \in J^+}\cl{L_k^{-1}(T_k-t)}{1} |\mr{\ups}_k (z)|^p\frac{c_j(X_{j,k}(z, T_k -t))}{c_k(z)}dz\right.\nn\\& \phantom{xxxxxxxxxxxxxx}+ \left.\sum\limits_{k \in J^-}\cl{0}{L_k^{-1}(t)}|\mr{\w}_k(z)|^p\frac{c_j(Y_{j,k}(z,t))}{c_k(z)}dz \right) \leq D_j \|(\mr {\mb\ups},\mr{\mb\w})\|^p_{\mb X_p}\nn,
     \end{align}
     where $D^{p/q}_b = \left(\sum\limits_{k \in J^+\cup J^-}(b_{jk})^q\right)^{p/q}$ and  $D_j$ is  a constant depending on $D^{p/q}_b$ and $\max_{k\in J^+\cup J^-}\left\{\frac{\max_{x\in [0,1]}c_j(x)}{\min_{x\in [0,1]}c_k(x)}\right\}$.

     In the same way we obtain
     \begin{equation}\label{west1}
        \cl{0}{1} |\w_j(x,t)|^pdx \leq   D_j \|(\mr {\mb\ups},\mr{\mb\w})\|^p_{\mb X_p}, \quad j\in J^-.
          \end{equation}
Hence, combining the estimates for all $j\in J^+\cup J^-$, we see that there is $D$ such that for all $(\mr{\mb\ups}, \mr{\mb \w})\in \mb X_p$ and $t\in [0,T]$,
    \begin{equation}
    \|\mc G(t)(\mr{\mb\ups}, \mr{\mb \w})\|_{\mb X_p} \leq D\|(\mr{\mb\ups}, \mr{\mb\w})\|_{\mb X_p}.\label{L2est1}
    \end{equation}

    Now, since $\mb X_p\subset \mb X_1$ and  $G_1(t)(\mr {\mb\ups},\mr{\mb\w}) = \mc G(t)(\mr {\mb\ups},\mr{\mb\w})$ for $(\mr {\mb\ups},\mr{\mb\w})\in D(A_{1,\mc B}),$ which is dense in both $\mb X_1$ and  $\mb X_p,$  we have
    $$
    G_1(t)|_{\mb X_p} = \mc G(t)|_{\mb X_p}, \quad 0\leq t\leq T,
    $$
     and $G_1(T)(\mr {\mb\ups},\mr{\mb\w}) = \mc G(T)(\mr {\mb\ups},\mr{\mb\w}).$ We can then  repeat the above procedure for $t = s+T, 0\leq s\leq T,$ getting $$\mc G(s)\mc G(T)(\mr {\mb\ups},\mr{\mb\w}) = G_1(s)G_1(T)(\mr {\mb\ups},\mr{\mb\w}) = G_1(t)(\mr {\mb\ups},\mr{\mb\w})\in \mb X_p,\; t= s+T, 0\leq s\leq T.$$
     Since \sem{G_1} is an algebraic semigroup, we see that, by iteration, $\sem{G_p}:=(G_1(t)|_{\mb X_p})_{t\geq 0}$ is an exponentially bounded (by \eqref{L2est1}) semigroup of continuous operators on $\mb X_p$. We have to prove its strong continuity at 0. For small $t>0$ we can use \eqref{upsdef} and \eqref{wdef}. Thus, for a given $j\in J^+,$
          \begin{align*}
     &\|\ups_j(\cdot,t)-\mr{\ups}_j\|^p_{L_p(0,1)} =\cl{L_j^{-1}(t)}{1}\left|\mr {\ups}_j(X_j(x,t)) - \mr\ups_j(x)\right|^p dx \\
     &+\!\!\!\!\cl{0}{L_j^{-1}(t)}\!\left|\sum\limits_{k \in J^+}\!\!b_{jk}\mr{\ups}_k (X_{k,j}(x,t-T_k))
  +\!\! \sum\limits_{k\in J^-}b_{jk}\mr{\w}_k(Y_{k,j}(x,t)) -\mr\ups_j(x)\right|^p \!\!dx\\
  &= I_1(t)+I_2(t).
     \end{align*}
     Consider first $\mr\ups_j \in C^\infty_0(0,1)$. Then, by \eqref{Lj}, $\mr\ups_j \circ L_j^{-1}$ is uniformly Lipschitz on $[0,L_j(1)]$ with a Lipschitz constant $\mu_j$. Since for $L_j^{-1}(t)\leq x\leq 1$ we have $0\leq L_j(x)-t \leq L_j(1)-t\leq L_j(1)$,  we can write
     $$
     I_1(t) \leq  \cl{L_j^{-1}(t)}{1}\left|\mr {\ups}_j(L_j^{-1}(L_j(x)-t)) - \mr\ups(L_j^{-1}(L_j(x)))\right|^p dx \leq \mu_j^p\cl{L_j^{-1}(t)}{1}t^p dx\leq \mu_j^p t^p.
     $$
      Then, as in \eqref{upse1},
          \begin{align*}
     &I_2(t) \leq  2^pD_b^{p/q} \left( \sum\limits_{k \in J^+}\cl{L_k^{-1}(T_k-t)}{1} |\mr{\ups}_k (z)|^p\frac{c_j(X_{j,k}(z, T_k-t))}{c_k(z)}dz\right.\\& + \left.\sum\limits_{k \in J^-}\cl{0}{L_k^{-1}(t)}|\mr{\w}_k(z)|^p\frac{c_j(Y_{j,k}(z,t))}{c_k(z)}dz \right) + 2^p\cl{0}{L_j^{-1}(t)}|\mr\ups_j(x)|^pdx \\
    & \leq 2^p D_j \left( \sum\limits_{k \in J^+}\cl{L_k^{-1}(T_k-t)}{1}\!\! |\mr{\ups}_k (z)|^pdz+ \sum\limits_{k \in J^-}\!\!\!\cl{0}{L_k^{-1}(t)}\!\!|\mr{\w}_k(z)|^pdz \right) + 2^p\!\!\!\!\cl{0}{L_j^{-1}(t)}\!\!\!|\mr\ups_j(x)|^pdx.
     \end{align*}
          Since $L_k^{-1}(t)\to 0$ and $L_k^{-1}(T_k-t)\to 1$ as $t\to 0^+$, $\lim_{t\to 0^+}I_2(t) =0$. The convergence can be extended to an arbitrary $\mr\ups_j \in L_p(0,1)$ by density and uniform boundedness, \eqref{L2est1}. The convergence of $\w_j$,  $j\in J^-,$ can be proved in the same way.
          Thus, $(\mr\ups, \mr \w) \in \mb X_p,$
    \begin{equation}
    \lim\limits_{t\to 0^+} G_p(t)(\mr{\mb\ups}, \mr{\mb \w}) = (\mr{\mb\ups}, \mr{\mb \w}),
    \end{equation}
    in $\mb X_p$ and hence \sem{G_p} is a strongly continuous semigroup in $\mb X_p$. Then the application of \cite[Proposition II.2.3]{EN} shows that the generator of \sem{G_p} is the part of $A_{1,\mc B}$ in $\mb X_p$. 
    \end{proof}

        \section{Examples}
        We shall discuss the relation of \eqref{bc2s} with the boundary conditions of \cite{Nic,KMN}.
    \begin{example} \label{ex38} Consider the model of \cite{Nic}
 \begin{equation}
  \p_tp^{j}_1 + K^j \p_x p^j_2  =0, \quad  \p_tp^{j}_2 + L^j \p_x p^j_1  =0,
      \label{sysN}
   \end{equation}
   for $t>0, 0<x<1, 0\leq j\leq m,$ where $K^j>0,L^j>0$  for all $j$.  In this case
   \begin{equation}
   \la^j_\pm = \pm\sqrt{L^jK^j}
   \label{speed}
   \end{equation}
   and we can set
   $$
   \mc F^j = \left(\begin{array}{cc}K^j&K^j\\\sqrt{L^jK^j}&-\sqrt{L^jK^j}\end{array}\right).
   $$
   For a given vertex $\bv,$ we introduce a function defined as $\nu^j(\bv) =  -1$ if $ l_j(\bv) = 0$ and $\nu^j(\bv) =  1$ if $l_j(\bv)= 1$
   and define $T_{\bv} \mb p_2 (\bv)= (\nu^j(\bv)p^j_2(\mb v))_{j\in J_{\bv}}.$
      In this case $\alpha_j=1$ for any $j$ and thus for any vertex $\bv$ we need $|J_{\bv}|$ boundary conditions. The ones introduced in \cite{Nic} can be expressed as follows. Let first $\bv$ be a vertex with $|E_{\bv}|>1$. We split $\mbb R^{|J_{\bv}|}$ into $X_{\bv}$ of dimension $n_{\bv}$ and its orthogonal complement $X_{\bv}^\perp$ of dimension $l_{\bv} = |J_{\bv}|-n_{\bv}$. Then it is required that
      $$
      \mb p_1(\bv) \in X_{\bv}, \quad T_v\mb p_2(\bv) \in X^\perp_{\bv},
      $$
that is, denoting $I_1=\{1,\ldots,n_{\bv}\}$ and $I_2=\{n_{\bv}+1,\ldots,|J_{\bv}|\}$,
\begin{equation}\label{Nbc}
\sum\limits_{j\in J_{\bv}}\phi^j_r p^j_1(\bv) = 0, \quad r \in I_2,\quad
\sum\limits_{j\in J_{\bv}}\vhi^j_r \nu^j(\bv)p^j_2(\bv) = 0, \qquad r \in I_1,
\end{equation}
where $((\vhi^j_r)_{j\in J_{\bv}})_{r\in I_1}$ is a base in $X_{\bv}$ and  $((\phi^j_r)_{j\in J_{\bv}})_{r\in I_2}$ is a base in $X^\perp_{\bv}$ so that $(\vhi^j_r)_{j\in J_{\bv}} \cdot (\phi^j_s)_{j\in J_{\bv}} = 0$ for any $r\in I_1$ and $s\in I_2$. In the notation of \eqref{Phiv}, we have $
\Phi^j_{\bv,r} = (\phi^j_r,0)$ for $r \in I_2,$ and $\Phi^j_{\bv,r}= (0,\nu^j(\bv)\vhi^j_r)$ for $r\in I_1$.
Thus
\begin{align*}
\Phi^j_{\bv,r}\cdot f^j_\pm &= \phi^j_r f^j_{\pm,1} = K^j\phi^j_r, \quad r \in I_2,\\
\Phi^j_{\bv,r} \cdot f^j_\pm &=\nu^j(\bv)\vhi^j_rf^j_{\pm,2} = -\sqrt{K^jL^j}\vhi^j_r, \quad r\in I_1.
\end{align*}
Then \eqref{bcsplit1} becomes
 \begin{equation}
 \begin{split}
   \sum\limits_{j \in J^0_{\bv}} \phi^j_rK^ju^j_1(0)
    +\sum\limits_{j\in J^1_{\bv}} \phi^j_r K^ju^j_2(1)&
     = \mc B_{in,r}, \quad r \in I_2,\\
    -\sum\limits_{j \in J^0_{\bv}} \vhi^j_r\sqrt{K^jL^j}u^j_1(0)
    -\sum\limits_{j\in J^1_{\bv}} \vhi^j_r\sqrt{K^jL^j}u^j_2(1)&= \mc B_{in,r}, \quad r \in I_1,
    \end{split}
    \label{bcsolv}
    \end{equation}
    where $\mc B_{in,r}$ are the incoming components in \eqref{bcsplit1}. Vectors $\{(\phi^j_rK^j)_{j\in J_{\bv}}\}_{r\in I_2}$  and $\{(-\vhi^j_r\sqrt{K^jL^j})_{j\in J_{\bv}}\}_{r\in I_1}$ form bases in, respectively, $X_{\bv}^\perp$ and $X_{\bv}$. Indeed, $\{(\phi^j_r)_{j\in J_{\bv}}\}_{r\in I_2}$ forms a basis in $X_{\bv}^\perp$ by assumption, thus the matrix $(\phi^j_r)_{{j\in J_{\bv}}, r\in I_2}$ has rank $l_{\bv}= |J_{\bv}|-n_{\bv},$ hence there is  a minor such that $\text{det} (\phi^{j_i}_r)_{{1\leq i\leq l_{\bv}}, r\in I_2}\neq 0.$ Then $\text{det} (\phi^{j_i}_rK^{j_i})_{{1\leq i\leq l_{\bv}}, r\in I_2}= \text{det} (\phi^{j_i}_r)_{{1\leq i\leq l_{\bv}}, r\in I_2}\prod\limits_{i=1}^{l_{\bv}}K^{j_i}\neq 0.$ The same argument is valid for $\{(-\vhi^j_r\sqrt{K^jL^j})_{j\in J_{\bv}}\}_{r\in I_1}$. Further,   $(\phi^j_rK^j)_{j\in J_{\bv}}$ is orthogonal to  $(-\vhi^j_s\sqrt{K^jL^j})_{j\in J_{\bv}}$ for any $r\in I_2$ and $s\in I_1$ with respect to the scalar product weighted with the vector $((K^j)^{-\frac{3}{2}}(L^j)^{-\frac{1}{2}})_{j\in J_{\bv}}$. Hence, the set $$\{\{(\phi^j_rK^j)_{j\in J_{\bv}}\}_{r\in I_2}, \{(-\vhi^j_r\sqrt{K^jL^j})_{j\in J_{\bv}}\}_{r\in I_1}\}$$
    forms a basis in $\mbb R^{|J_{\bv}|}$ and therefore \eqref{bcsolv} is uniquely solvable.

    To complete the analysis, let us consider the vertices with $|E_{\bv}|=1,$ referred to  as exterior in \cite{Nic}.  The author splits arbitrarily the set of exterior vertices into $V_{ext}^{\mathrm{Diss}}$ and $V_{ext}^{\mathrm{Dir}},$ whereupon we require
    $$
    p_2^{j_{\bv}} (\bv) = 0, \quad \bv \in V_{ext}^{\mathrm{Dir}}, \qquad   p_1^{j_{\bv}}(\bv,t) = \alpha \nu^{j_{\bv}}(\bv)p_2^{j_{\bv}} (\bv),\quad \bv \in V_{ext}^{\mathrm{Diss}},
    $$
    where $\alpha\geq 0$. The first condition can be re-written as
    $$
    f^{j_{\bv}}_{+,2}u_1^{j_{\bv}}(\bv) + f^{j_{\bv}}_{-,2}u_2^{j_{\bv}}(\bv) = 0.
    $$
    This corresponds to $\Phi^{j_{\bv}}_{\bv} = (0,1)$ and, since $f^{j_{\bv}}_{+,2} = \sqrt{K^{j_{\bv}}L^{j_{\bv}}} = - f^{j_{\bv}}_{-,2} \neq 0$, the equation above is solvable for either $u_1^{j_{\bv}}(0)$ or $u_2^{j_{\bv}}(1),$ whichever is necessary. The second condition can be re-written as
$$
(K^{j_{\bv}} -\alpha\nu^{j_{\bv}}(\bv) \sqrt{K^{j_{\bv}}L^{j_{\bv}}})u_1^{j_{\bv}}(\bv) + (K^{j_{\bv}} +\alpha \nu^{j_{\bv}}(\bv)\sqrt{K^{j_{\bv}}L^{j_{\bv}}})u_2^{j_{\bv}}(\bv) =0.
$$
Thus, if $l_{j_{\bv}}(\bv)=0$, then $\nu^{j_{\bv}}(\bv) =-1$ and $(K^{j_{\bv}} -\alpha\nu^{j_{\bv}}(\bv) \sqrt{K^{j_{\bv}}L^{j_{\bv}}})\neq 0,$ yielding the solvability of the equation with respect to   $u_1^{j_{\bv}}(\bv).$ Similarly, if $l_{j_{\bv}}(\bv)=1$, then $\nu^{j_{\bv}}(\bv) =1$ and the equation is solvable for $u_2^{j_{\bv}}(\bv)$.

To conclude this example, let us consider the approach introduced in \cite{KMN}, where the boundary conditions are imposed after the matrix $\mc M^j$ has been transformed to a Hermitian form via a Hermitian matrix $Q^j$. For a given vertex $\bv$ the authors consider the block-diagonal matrix $
T_{\bv} := (Q^{j}\mc M^{j}\iota^{j} (\bv))_{j \in J_{\bv}},
$
where $\iota^{j} (\bv)$ is the diagonal matrix $\textrm{diag}\{\nu^{j}(\bv)\}$. Then the boundary conditions for a function $\Theta= (\theta^j)_{j\in E_{\bv}},$ where $\theta^j$ is a vector of functions of dimension equal to the number of equations on $\mb e^{j}$,  are formulated as the requirement that $\Theta(\bv)$ belongs to the totally isotropic subspace associated with the quadratic form $T_{\bv}\Theta(\bv )\cdot \bar \Theta(\bv).
$
In the considered case a symmetrizing matrix $Q^j$ and the symmetrization are given by, respectively, $$
Q^{j} =\left(\begin{array}{cc} L^{j}&0\\0&K^{j}\end{array}\right),\qquad
Q^{j} \mc M^{j} = \left(\begin{array}{cc} 0&L^{j}K^{j}\\L^{j}K^{j}&0\end{array}\right).
$$
Restricting our attention to real solutions and denoting $\theta^{j} = (p_1^{j},p_2^{j})$, the boundary condition can be written as
\begin{equation}
 -\sum\limits_{j\in J_{\bv}^0} K^jL^j p^j_1(0)p^j_2(0) + \sum\limits_{j\in J_{\bv}^1} K^jL^j p^j_1(1)p^j_2(1) = 0.
\label{KMN}
\end{equation}
If we consider an example of the boundary conditions discussed in \cite{Nic, KMN}
\begin{equation}
p_1^{j}\;\text{are\; continuous\; across\;} \bv\; \text{for}\; j \in J_{\bv}\; \text{and}\; \sum\limits_{j\in J_{\bv}} \nu^{j}(\bv) p_2^{j}(\bv) =0,
\label{nic}
\end{equation}
then we need $p^j_1(\bv) = p$ for some $p\in \mbb R$ and all $j\in J_{\bv}$ as well as
\begin{equation}
-\sum\limits_{j\in J_{\bv}^0}  p_2^{j}(0) + \sum\limits_{j\in J_{\bv}^1}  p_2^{j}(1)=0.\label{KMN1}
\end{equation}
On the other hand, after dividing by $p$, \eqref{KMN} requires
\begin{equation}
-\sum\limits_{j\in J_{\bv}^0} K^jL^jp^j_2(0) + \sum\limits_{j\in J_{\bv}^1} K^jL^j p^j_2(1) = 0.
\label{KMN2}
\end{equation}
Thus, the subspaces determined by \eqref{KMN1} and \eqref{KMN2} coincide only if $K^jL^j = KL$ for all $j \in J_{\bv}$, that is, the boundary conditions of \cite{KMN} cover the Kirchhoff law only when the speeds $\sqrt{L^jK^j},$ see \eqref{speed}, are the same on all edges.
\end{example}
\begin{example} \label{KMNex}
To further compare the approach of \cite{KMN} and of this paper, consider a simple one edge network case with the telegraph equation
\begin{equation}\label{exte}
\p_t p_1 = \p_x p_2, \quad \p_t p_2 = \p_x p_1
\end{equation}
in $(L_2(0,1))^2$. The edge $\mb e$ is identified with $(0,1)$. Let $\mb v_1$ be the endpoint associated with 0, and $\mb v_2$ the endpoint associated with 1.
Then
\begin{equation}
\begin{split}
\Re \int_0^1 (\p_x p_2,\p_x p_1)\cdot {(p_1,p_2)}dx &= \Re( p_1(1)\bar p_2(1)) - \Re (p_1(0)\bar p_2(0))\\& =: q_{\mb v_2}(\mb p(\mb v_2)) + q_{\mb v_1}(\mb p(\mb v_1)),\end{split}
\label{diss}
\end{equation}
where $q_{\mb v_i}, i=0,1,$ are quadratic forms.  The boundary conditions considered in \cite{KMN} are expressed as
$(p_1(\mb v), p_2(\mb v))\in Y_{\mb v}\subset \mbb C^2, \mb v= \mb v_1,\mb v_2,$ where $Y_{\mb v}$ is a linear space on which $q_{\mb v}$ vanishes, called a total isotropic subspace associated with $q_{\mb v}$. Equivalently, the boundary conditions can be expressed as
\begin{equation}
\mb w^{(\mb v,i)}\cdot {(p_1(\mb v), p_2(\mb v))} = 0,
\label{yperp}
\end{equation}
where $\{\mb w^{v,i}\}_{1\leq i\leq dim Y^\perp_{\mb v}}$ is a fixed basis in $Y^\perp_{\mb v}\subset \mbb C^2;$ it is possible that $Y^\perp_{\mb v}=\{0\}$ in which case the above system does not impose any conditions at $\mb v$.  The basic solvability condition of \cite[Lemma 3.5]{KMN}, specified to this case, is
\begin{equation}
dim Y_{\mb v_1} + dim Y_{\mb v_2} = dim (Y_{\mb v_1}^\perp + Y_{\mb v_2}^\perp) =2,
\label{KMNcond}
\end{equation}
where the sum in the middle is the algebraic sum of the subspaces and both $Y_{\mb v_1}^\perp$ and $Y_{\mb v_2}^\perp$ are considered as subspaces of $\mbb C^2$.

If we consider the boundary conditions $p_1(0) = p_1(1) =0$ (with which \eqref{exte} reduces to the wave equation for $p_1$ with the homogeneous Dirichlet boundary conditions, see \cite{Ban97}), then $Y_{\mb v_i} = \mc {Lin}\{(0,1)\}, i=1,2,$ (where $\mc{Lin}$ denotes the linear span) and thus $dim Y_{\mb v_1} + dim Y_{\mb v_2}=2.$ However,  $Y_{\mb v_1}^\perp = Y_{\mb v_2}^\perp =\mc {Lin}\{(1,0)\}$ and $dim (Y_{\mb v_1}^\perp + Y_{\mb v_2}^\perp) = 1$. On the other hand, if we consider the boundary conditions $p_1(0) = p_2(1) =0$ (with which \eqref{exte} reduces to the wave equation for $p_1$ with the Dirichlet boundary condition at $x=0$ and the Neuman condition at $x=1$), then $ Y_{\mb v_1} = \mc {Lin}\{(0,1)\},  Y_{\mb v_2} = \mc {Lin}\{(1,0)\},$ and thus again $dim Y_{\mb v_1} + dim Y_{\mb v_2}=2.$ However,  $Y_{\mb v_1}^\perp = \mc {Lin}\{(1,0)\},$ while $Y_{\mb v_2}^\perp =\mc {Lin}\{(0,1)\}$ and $dim (Y_{\mb v_1}^\perp + Y_{\mb v_2}^\perp) = 2$.

Due to \eqref{diss}, both boundary conditions lead to the problem being dissipative but the first one does not satisfy the solvability condition \eqref{KMNcond}.   This rather restrictive condition allows for checking $m$-dissipativity of the problem by simply considering the resolvent equation with $\la =0$,
\begin{equation}
\p_x p_2 = f_1, \quad \p_x p_1 = f_2.
\label{0max}
\end{equation}
Indeed, since the resolvent set is open, if it contains $\la =0$, then it also contains some $\la>0$ and thus the assumptions of the Lumer-Phillips theorem, \cite[Theorem II.3.15]{EN}. This significantly  simplifies the analysis in the more general network case. Indeed, here
$$
p_2(x) = K^1 + \int_0^x f_1(s)ds, \quad p_1(x) = K^2 + \int_0^x f_2(s)ds
$$
and for $p_1(0)=p_1(1)=0$, conditions \eqref{yperp} take the form  of the system
\begin{align*}
(1,0)\cdot {(p_1(0),p_2(0))} &= (1,0)\cdot ( K^2, K^1) =0,\\
(1,0)\cdot {(p_1(1),p_2(1))} &= (1,0)\cdot ( K^2, K^1) +(1,0)\cdot \left(\int_0^1  f_2(s)ds,\int_0^1  f_1(s)ds\right)=0,
\end{align*}
not solvable for general $(f_1,f_2)\in (L_2(0,1))^2$. If $p_1(0)=p_2(1)=0$, we get
\begin{align*}
(1,0)\cdot {(p_1(0),p_2(0))} &= (1,0)\cdot ( K^2, K^1) =0,\\
(0,1)\cdot (p_1(1),p_2(1)) &= (0,1)\cdot ( K^2, K^1) +(0,1)\cdot \left(\int_0^1  f_2(s)ds,\int_0^1  f_1(s)ds\right)=0
\end{align*}
which has a solution.

If we consider the full resolvent equation
\begin{equation}
\la p_1-\p_x p_2 = f_1, \quad \la p_2-\p_x p_1 = f_2,
\label{lmax}
\end{equation}
then, adding and subtracting the equations to diagonalize the system, we obtain the general solution as
\begin{align*}
p_1(x)&=   \frac{A^1e^{\la x}+ A^2e^{-\la x}}{2}\\
& - \frac{ \int_0^x (e^{\la (x-s)}-e^{-\la(x-s)})f_1(s)ds +\int_0^x (e^{\la (x-s)}+e^{-\la(x-s)})f_2(s))ds}{2},\\
p_2(x)&=   \frac{A^1e^{\la x}- A^2e^{-\la x}}{2}\\
& - \frac{ \int_0^x (e^{\la (x-s)}+e^{-\la(x-s)})f_1(s)ds +\int_0^x (e^{\la (x-s)}-e^{-\la(x-s)})f_2(s))ds}{2} .
\end{align*}
It is easy to see that if we impose $p_1(0)=p_1(1) =0,$ determining $A^1$ and $A^2$ requires invertibility of $\left(\begin{array}{cc}1&1\\e^\la&e^{-\la}\end{array}\right).$ Clearly, $\la =0$ is an eigenvalue  (which explains the failure of previous approach) but any $\la \in \mbb R\setminus \{0\}$ is in the resolvent set and hence we have a generation of a group. Similarly, imposing $p_1(0)=p_2(1) =0$  requires invertibility of $\left(\begin{array}{cc}1&1\\e^\la&-e^{-\la}\end{array}\right)$ and in this case $\la =0$ is not an eigenvalue and the $m$-dissipativity can be determined from \eqref{0max}.

We observe that $u_1 = p_1+p_2$ and $u_2=p_1-p_2$ are the Riemann invariants of the system and satisfy
\begin{equation}\label{exte1}
\p_t u_1 - \p_x u_1 =0, \quad \p_t u_2 +\p_x u_2 =0
\end{equation}
on $(0,1)$. Clearly, $u_2$ flows from 0 to 1 and $u_1$ from 1 to 0. Hence, $u_2(0)$ is outgoing at $x=0$, while $u_1(1)$ is outgoing at $x=1$ and the conditions $p_1(0)=p_1(1)=0$ in terms of matrices $\mb \Psi_{\mb v_1}, \mb \Psi_{\mb v_2}$ can be written as
$$
\mb \Psi_{\mb v_1}\binom {u_1(0)}{u_2(0)} = (1,1)\binom {u_1(0)}{u_2(0)} = 0, \quad \mb \Psi_{\mb v_2}\binom {u_1(1)}{u_2(1)} = (1,1)\binom {u_1(1)}{u_2(1)} = 0
$$
and assumptions of Proposition \ref{prop1} are satisfied.
\end{example}
\begin{example}\label{ex39}
Let us consider the linearized Saint-Venant system \eqref{SV3l} (where we replaced the original variables by $(p^j_1,p^j_2)$),
\begin{equation}
\p_tp^j_1 = -V^j \p_x p^j_1 - H^j\p_x p^j_2, \quad \p_t p^j_2= -g \p_x p^j_1 -V^j \p_x p^j_2,
\label{SV}
\end{equation}
and assume that on each edge we have $Fr^j >1$; that is, $\la^j_{\pm} = V^j \pm \sqrt{gH^j} >0$. If we consider the network shown on Fig. \ref{fig1}, we see that at $\bv_0$ we need two boundary conditions, no boundary conditions at $\bv_i$, $i=2,\ldots,N$ and $2N-2$ boundary conditions at $\bv_1$.
In this case we have
\begin{equation}
 \binom{p^j_1}{p^j_2} = \binom{f^j_{+,1} u^j_1 + f^j_{-,1}u^j_2}{f^j_{+,2}u^j_1 + f^j_{-,2}u^j_2} = \binom{H^j u^j_1+ H^j u^j_2}{\sqrt{gH^j}u^j_1 -\sqrt{gH^j}u^j_2}.\label{SNpguw}
 \end{equation}

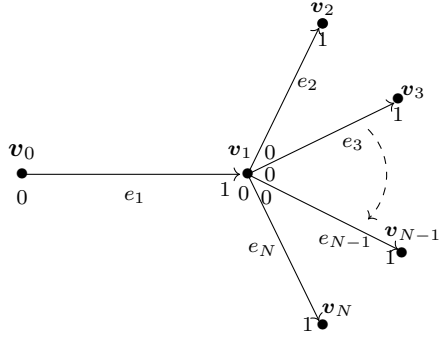
\begin{figure}
\centering
\begin{tikzpicture}
\draw [->] (-3,0) -- (-0.1,0);
\draw  [dashed, ->] (1.6,0.6) to [out=-45,in=45] (1.6,-0.6);
\draw (-0.1,0.3) node {\footnotesize{$\bv_1$}};
\draw(0,0)node{$\bullet$};
\draw (-3,0)node{$\bullet$};
\draw  [->](0,0)--(1.95,0.95);
\draw  [->](0,0)--(2,-1);
\draw[->](0,0)--(0.95,1.95);
\draw [->](0,0)--(0.95,-1.95);
 \draw (-3,-0.3) node {\footnotesize{$0$}};
 \draw (-3,0.3) node {$\bv_0$};
 \draw (1.4,0.4) node {\footnotesize{$e_3$}};
 \draw (1.3,-0.9) node {\footnotesize{$e_{N-1}$}};
 \draw (-1.5,-0.3) node {\footnotesize{$e_1$}};
 \draw (-0.3,-0.2) node {\footnotesize{1}};
  \draw (0.3,0.3) node {\footnotesize{$0$}};
  \draw (0.3,0) node {\footnotesize{$0$}};
  \draw (0.25,-0.3) node {\footnotesize{$0$}};
  \draw (-0.05,-0.25) node {\footnotesize{$0$}};
 \draw (1,2.2) node {\footnotesize{$\bv_2$}};
 \draw (1,1.8) node {\footnotesize{$1$}};
 \draw(1,2)node{$\bullet$};
 \draw(2,1)node{$\bullet$};
 \draw(2.05,-1.05)node{$\bullet$};
 \draw(1,2)node{$\bullet$};
 \draw (1.2,-1.8) node {\footnotesize{$\bv_{N}$}};
 \draw (2.2,-0.8) node {\footnotesize{$\bv_{N-1}$}};
 \draw (2.2,1.1) node {\footnotesize{$\bv_3$}};
 \draw(1,-2)node{$\bullet$};
 \draw (0.8,1.2) node {\footnotesize{$e_2$}};
 \draw (0.2,-1) node {\footnotesize{$e_{N}$}};
 \draw (1.9,-1.1) node {\footnotesize{$1$}};
 \draw (2,0.8) node {\footnotesize{$1$}};
 \draw (0.8,-2) node {\footnotesize{$1$}};
  \end{tikzpicture}
\caption{Starlike network of channels}\label{fig1}
\end{figure}
Here $J_{\bv_0}=J_{\bv_0,2} = J^0_{\bv_0} =\{1\}$ and all other subsets of $J_{\bv_0}$ are empty. Hence \eqref{bc2uw} takes the form
    $$
    \left(\begin{array}{cc}\phi^1_{\bv_0,1}&\vhi^1_{\bv_0,1}\\\phi^1_{\bv_0,2}&\vhi^1_{\bv_0,2}\end{array}\right)\left(\begin{array}{cc} f^j_{+,1}&f^j_{-,1}\\
   f^j_{+,2}&f^j_{-,2} \end{array}\right)\binom{u^1_1(0)}{u^1_2(0)} = \binom{0}{0}
$$
and the equation is solvable if and only if $\Phi^1_{\bv_0,1}$ and $\Phi^1_{\bv_0,2}$ are linearly independent, on account of the invertibility of $\mc F^1$. Similarly, for $\bv_1,$ \eqref{bc2} takes the form
\begin{equation}
  \sum\limits_{j =1}^{N} (\phi^j_{\bv_1,r} p^j_1(\bv_1) + \vhi^j_{\bv_1,r} p^j_2(\bv_1))=0, \quad r=1,\ldots, 2N-2
  \label{bc3}
  \end{equation}
and \eqref{compsolv} is satisfied if and only if
$$
\widetilde{\mb\Phi_{\bv_1}} \text{diag}\{\mc F^j\}_{2\leq j\leq N},
$$
where
\begin{align*}
\widetilde{\mb \Phi_{\bv_1}}:& = 
\left(\begin{array}{ccccc}\phi^2_{\bv_1,1}&\vhi^2_{\bv_1,1}&\ldots&\phi^N_{\bv_1,1}&\vhi^N_{\bv_1,1}\\
\vdots&\vdots&\vdots&\vdots&\vdots\\
\phi^2_{\bv_1,2N-2}&\vhi^2_{\bv_1,2N-2}&\ldots&\phi^N_{\bv_1,2N-2}&\vhi^N_{\bv_1,2N-2}\end{array}\right),
\end{align*}
is nonsingular; that is,  if and only if $\widetilde{\mb \Phi_{\bv_1}}$
is invertible. In the case of \cite{KMN},
$$
p^j_1(0) = p^1_1(1), \quad p^j_2(0) = p^1_2(1),\qquad  j=2,\ldots, N,
$$
we have
\begin{align*}
\Phi_{\mb v_1,r} &= (-1,0,0,0,\ldots,1,0,0,0,\ldots,0,0),\quad r=1,3,\ldots 2N-3,\\
\Phi_{\mb v_1,r} &= (0,-1,0,0,\ldots,0,1,0,0,\ldots,0,0),\quad r=2,4,\ldots 2N-2,
\end{align*}
where in both cases $1$ appears at the $(r+2)$-th place. In this case $\widetilde{\mb\Phi_{\bv_1}} = I$ and thus the system of boundary conditions can be solved for $(u^j_1(0),u^j_2(0))_{2\leq j\leq N}$.
\end{example}

\bibliographystyle{abbrv}
\bibliography{research}

 \end{document}